\documentclass[11pt,a4paper]{article}

\usepackage{amsfonts,amsmath,amssymb,amsopn}
\usepackage{vmargin}

\usepackage[all]{xy}

\newtheorem{thm}{Theorem}[section]
\newtheorem{prop}[thm] {Proposition}
\newtheorem{cor}[thm] {Corollary}
\newtheorem{lem}[thm] {Lemma}

\newtheorem{dfn}[thm] {Definition}

\title{A Simpson correspondence in positive characteristic}
\author{Michel Gros, Bernard Le Stum \& Adolfo Quir\'os\footnote{A.~Q. was partially supported by Project GALAR (MTM2006-10548)
from MEC (Spain) and by the
joint Madrid Region-UAM project TENU2 (CCG07-UAM/ESP-1814).}}
                                    
\begin{document}
\maketitle

\tableofcontents

\section*{Introduction}

The first two authors had the opportunity to participate in a working group in Rennes dedicated to the work of Arthur Ogus and Vadim Vologodsky on non abelian Hodge theory, which is now published in  \cite{OgusVologodsky07}. This is an analog in positive characteristic $p$ of Simpson's correspondence over the complex numbers between local systems and a certain type of holomorphic vector bundles that he called Higgs bundles (\cite{Simpson92}).
Actually, Pierre Berthelot had previous results related to these questions and he used this opportunity to explain them to us.
What we want to do here is to extend these results to differential operators of higher level.
In the future, we wish to lift the theory modulo some power of $p$ and compare to Faltings-Simpson $p$-adic correspondence (\cite{Faltings05}).

We would also like to mention some other papers related to our investigation.
First, there is an article \cite{Kaneda04} by Masaharu Kaneda where he proves the (semi-linear) Azumaya nature of the ring of differential operators of higher level, generalizing the result of \cite{BezrukavnikovMirkovicRumynin08}.
Also, Marius van der Put in \cite{Put95} studies the (linear) Azumaya nature of differential operators in the context of differential fields.

We will first recall, in sections \ref{usual} and \ref{higher}, in an informal way the notion of divided powers of higher level and how this leads, using some duality, to arithmetic differential operators.
Then, in section \ref{pem}, we will define the notion of $p^m$-curvature and show that Kaneda's isomorphism still holds over an arbitrary basis.
Next, in section \ref{lift}, we assume that there exists a strong lifting of Frobenius mod $p^2$ and use it to lift the divided Frobenius and derive a Frobenius map on the ring of differential operators of level $m$.
Actually, we can do better and prove in Theorem \ref{mainth} that this data determines a splitting of a central completion of this ring.
It is then rather formal, in section \ref{higgs}, to obtain a Simpson correspondence and we can even give some explicit formulas.
We finish the article, in section \ref{comp},  with a series of complements concerning compatibility with other theories.

\section*{Acknowledgments}

Many thanks to Pierre Berthelot who explained us the level zero case and helped us understand some tricky constructions.

\section*{Conventions}

We let $p$ be a prime and $m \in \mathbf N$.
Actually, we are interested in the $m$-th power of $p$.
When $m = 0$, we have $p^m = 1$ which is therefore independent of $p$.
We may also consider the case $m = \infty$ in which case we will write $p^m = 0$.
Again, this is independent of $p$.

Unless $m = 0$ or $m = \infty$, all schemes are assumed to be $\mathbf Z_{(p)}$-schemes.

We use standard multiindex notations, hoping that everything will be clear from the context.

\section{Usual divided powers}\label{usual}

It seems useful to briefly recall here some basic results on usual divided powers that we will need afterwards.
There is nothing new but we hope that it makes the next sections easier to read without referring to older articles.
The main point in this section is to clarify the duality between divided powers and regular powers.

Let $\mathcal R$ be a commutative ring (in a topos).

\begin{dfn}
A \emph{divided power structure} on an ideal $\mathcal I$ in a commutative $\mathcal R$-algebra $\mathcal A$ is a family of maps
$$
\xymatrix@R=0cm{\mathcal I \ar[r] & \mathcal A \\ \quad f \ar@{|->}[r] &  f^{[k]}}
$$
that behave like $f \mapsto \frac {f^k}{k!}$.
We then say that $\mathcal I$ is a \emph{divided power ideal} or that $\mathcal A$ is a \emph{divided power $\mathcal R$-algebra}.
\end{dfn}

We will not list all the required properties.
Note however that we always have
$$
f^{[k]} f^{[l]} = {k+l \choose k} f^{[k+l]}.
$$

\begin{thm}
The functor $\mathcal A \mapsto \mathcal I$ from divided power $\mathcal R$-algebras to $\mathcal R$-modules has a left adjoint $\mathcal M \mapsto \Gamma_{\bullet} \mathcal M$.
\end{thm}

\textbf{Proof :} See for example, Theorem 3.9 of \cite{BerthelotOgus78}. $\Box$

Actually $\Gamma_{\bullet} \mathcal M$ is a graded algebra with divided power ideal $\Gamma_{>0} \mathcal M$.
Not also that  $\Gamma_{\bullet} \mathcal M$ is generated as $\mathcal R$-algebra by all the $s^{[k]}$ for $s \in M$.
For example, if $\mathcal M$ is free on $\{s_{\lambda}, \lambda \in \Lambda\}$, then $\Gamma_{k} \mathcal M$ is free on the $\underline s ^{[\underline k]} : = \prod s_{\lambda}^{[k_{\lambda}]}$ with $|\underline k| := \sum k_{\lambda}= k$.
Moreover, multiplication is given by the general formula for divided powers recalled above.

\begin{dfn} If $\mathcal M$ is an $\mathcal R$-module, the ring $\Gamma_{\bullet} \mathcal M$ is called the \emph{divided power algebra} on $\mathcal M$.
\end{dfn}

If $\mathcal M$ is an $\mathcal R$-module, we will denote by $S^{\bullet} \mathcal M$ the \emph{symmetric algebra} on $\mathcal M$ and by $\check {\mathcal M}$ the dual of $\mathcal M$.
Also, we will denote by $\widehat{S^{\bullet} \mathcal M}$ the completion of  $S^{\bullet} \mathcal M$ along $S^{>0} \mathcal M$.

\begin{prop}
If $\mathcal M$ is an $\mathcal R$-module, there exists a canonical pairing
$$
\xymatrix@R=0cm{S^{\bullet} \check{\mathcal M} \times \Gamma_\bullet \mathcal M \ar[r] &\mathcal R \\ (\varphi_{1}\cdots \varphi_{n}, s^{[n]}) \ar@{|->}[r] & \varphi_{1}(s)\cdots \varphi_{n}(s)}
$$
giving rise to perfect duality at each step when $\mathcal M$ is locally free of finite type.
\end{prop}

\textbf{Proof :} See for example, Proposition A.10 of \cite{BerthelotOgus78}). $\Box$

The general formula for this pairing is quite involved but if $\{s_{\lambda}, \lambda \in \Lambda\}$ is a finite basis for $\mathcal M$, and $\{\check s_{\lambda}, \lambda \in \Lambda\}$ denotes the dual basis, then the dual basis to $\{\underline {\check s}^{\underline k}\}$ is nothing else but $\{\underline s ^{[\underline k]}\}$.

\begin{cor}
If $\mathcal M$ a locally free $\mathcal R$-module of finite type, we have a perfect pairing
$$
\xymatrix@R=0cm{\widehat {S^{\bullet} \check{\mathcal M}} \times \Gamma_\bullet \mathcal M \ar[r] &\mathcal R}.
$$
\end{cor}

Of course, there exists also a natural map $S^{\bullet} \mathcal M \to \Gamma_\bullet \mathcal M$ but it is not injective in general: if $p^{N+1} = 0$ on $X$, then $f^{p^n} \mapsto p^nf^{[p^n]} = 0$ for $n > N$.

\begin{prop} \label{dualmul} If $\mathcal M$ is an $\mathcal R$-module, multiplication on $S^{\bullet} \check{\mathcal M}$ is dual to the diagonal map
$$
\xymatrix@R=0cm{\Gamma_\bullet \mathcal M \ar[r]^-\delta &  \Gamma_\bullet (\mathcal M \oplus \mathcal M) \ar[r]^-\simeq & \Gamma_\bullet \mathcal M \otimes \Gamma_\bullet \mathcal M \\ s^{[k]} \ar@{|->}[rr] && \sum_{i+j=k} s^{[i]} \otimes s^{[j]}}
$$
\end{prop}

\textbf{Proof :} The point is to show that $(\varphi_{1}\cdots \varphi_{i} \otimes \psi_{1} \cdots \psi_{j}) \circ \delta$ acts like
$\varphi_{1}\cdots \varphi_{i}\psi_{1} \cdots \psi_{j}$ on $s^{[k]}$ when $k = i + j$.
And this is clear.
$\Box$

\section{Higher divided powers}\label{higher}

We quickly recall the definition of  Berthelot's divided powers of level $m$ and how one derives the notion of differential operators of higher level from them (see \cite{Berthelot96} for a detailed exposition).
We stick to a geometric situation.

\begin{dfn}
Let $X \hookrightarrow Y$ be an immersion of schemes defined by an ideal $\mathcal I$.
A \emph{divided power structure of level $m$} on $\mathcal I$ is a divided power ideal $\mathcal J \subset \mathcal O_{Y}$ such that
$$
\mathcal I^{(p^m)} + p\mathcal I \subset \mathcal J \subset \mathcal I.
$$
\end{dfn}
Here, $\mathcal I^{(p^m)}$ denotes the ideal generated by $p^m$-th powers of elements of $\mathcal I$.

It is then possible to define \emph{partial divided powers} on $\mathcal I$ : they are maps
$$
\xymatrix@R=0cm{\mathcal I \ar[r] & \mathcal A \\ \quad f \ar@{|->}[r] &  f^{\{k\}}}
$$
that behave like $f \mapsto \frac {f^k}{q!}$ where $q$ is the integral part of $\frac k{p^m}$.
Actually, if $k = qp^m+r$ and $f \in \mathcal I_{X}$, one sets
$$
f^{\{k\}} := f^r (f^{p^{m}})^{[q]}.
$$
We have, as above, a multiplication formula (writing $q_{k}$ instead of $q$ in order to take into account the dependence on $k$):
$$
f^{\{k\}} f^{\{l\}} = \left\{{k+l \atop k}\right\} f^{\{k+l\}} \quad \mathrm{where} \quad \left\{{k+l \atop k}\right\} = \frac {q_{k+l}!}{q_{k}! q_{l}!}.
$$

When $m = 0$, we must have $\mathcal J = \mathcal I$ and $f^{\{k\}} = f^{[k]}$ is the above divided power.
When $m = \infty$, the condition reduces to $p\mathcal J \subset \mathcal I \subset \mathcal J$ and we may always choose $\mathcal J = \mathcal I$ since $p\mathcal I$ has divided powers.
Also, in this case, $f^{\{k\}} = f^k$ is just the usual power.

We fix a (formal) scheme $S$ with a divided power structure of level $m$ on some ideal of $\mathcal O_{S}$ and we assume that all constructions below are made over $S$ and are ``compatible'' with the divided powers on $S$ in a sense that we do not want to make precise here (see \cite{Berthelot96} for details).
Actually, in the case of a regular immersion, the divided power envelope defined below does not depend on $S$ and this applies in particular to the diagonal embedding of a smooth $S$-scheme.

\begin{prop}
If an immersion $X \hookrightarrow Y$ has divided powers of level $m$, there is a finest (decreasing) ring filtration $\mathcal I^{\{n\}}$ such that  $\mathcal I^{\{1\}} = \mathcal I$ and such that $f^{\{h\}} \in \mathcal I^{\{nh\}}$ whenever $f \in \mathcal I^{\{n\}}$.
\end{prop}

\textbf{Proof :} See Proposition 1.3.7 of \cite{Berthelot96}. $\Box$

\begin{prop}
The functor that forgets the divided power structure on an immersion $X \hookrightarrow Y$ has a left adjoint $X \hookrightarrow P_{Xm}(Y)$.
\end{prop}
\textbf{Proof :} See proposition 1.4.1 of \cite{Berthelot96}. $\Box$

\begin{dfn}
If $X \hookrightarrow Y$ is an immersion of schemes, then $P_{Xm}(Y)$ is called the \emph{divided power envelope} of level $m$ of $X$ in $Y$.
\end{dfn}

We will denote by $\mathcal P_{Xm}(Y)$ the structural sheaf of $P_{Xm}(Y)$ and by $\mathcal I_{Xm}(Y)$ the divided power ideal of level $m$.
We will also need to consider the usual divided power ideal $\mathcal J_{Xm}(Y)$ (in \cite{Berthelot96}, these ideals are denoted by $\overline {\mathcal I}$ and $\widetilde {\mathcal I}$ if $\mathcal I$ denotes the ideal that defines the immersion).
For each $n$, we will denote by $P^n_{Xm}(Y)$ the subscheme defined by $\mathcal I^{\{n+1\}}_{Xm}(Y)$ and consider its  structural sheaf
$$
\mathcal P^n_{Xm}(Y) = \mathcal P_{Xm}(Y)/\mathcal I^{\{n+1\}}_{Xm}(Y).
$$

We will mainly be concerned with diagonal immersions $X \hookrightarrow X \times_{S} X$, and we will then write $P_{Xm}$, $\mathcal P_{Xm}$, $\mathcal I_{Xm}$, $\mathcal J_{Xm}$,$P^{n}_{Xm}$ and $\mathcal P^{n}_{Xm}$ respectively.

If we are given local coordinates $t_{1}, \ldots, t_{r}$ on $X/S$, the ideal $\mathcal I$ of the diagonal immersion is generated by the $\tau_{i} = 1 \otimes t_{i} - t_{i} \otimes 1$.
We always implicitly use the first projection as structural map and therefore write $t_{i} \otimes 1 = t_{i}$ and $1 \otimes t_{i} = t_{i} + \tau_{i}$.
When $m = 0$, $\mathcal P_{Xm}$ is nothing but the divided power algebra on the free $\mathcal O_{X}$-module on the generators $\tau_{1}, \ldots, \tau_{r}$.
Of course, for $m = \infty$, this is just the symmetric algebra.

In general, we obtain
$$
\mathcal O_{X} \langle \tau_{1}, \ldots, \tau_{r} \rangle ^{(m)} := \{\sum_{\mathrm{finite}} f_{\underline i} \mathcal \underline \tau^{\{\underline i\}}, \quad f_{\underline i} \in \mathcal O_{X}\}
$$
with multiplication given by the general formula for divided powers of level $m$ recalled above.

\begin{dfn}
If $X$ is an $S$-scheme, the dual to $\mathcal P^{n}_{Xm}$ is the sheaf $\mathcal D^{(m)}_{X,n}$ of \emph{differential operators of level $m$ and order at most $n$} and $\mathcal D^{(m)}_{X} = \cup_{n} \mathcal D^{(m)}_{X,n}$ is the sheaf of differential operators of level $m$ on $X/S$.
\end{dfn}

There is a composition law on $\mathcal D^{(m)}_{X}$ that comes by duality from the morphism
$$
\xymatrix@R=0cm{\mathcal P_{Xm} \ar[r]^-\delta & \mathcal P_{Xm} \otimes \mathcal P_{Xm} \\ a \otimes b \ar[r] & ( a \otimes 1) \otimes (1 \otimes b).}
$$
When $X/S$ is smooth, this turns $\mathcal D^{(m)}_{X}$ into a non commutative ring.

Locally, we see that
$$
\mathcal D^{(m)}_{X} = \{\sum_{\mathrm{finite}} f_{\underline i} \mathcal \underline \partial^{<\underline i>}, \quad  f_{\underline i} \in \mathcal O_{X}\}
$$
where $\underline \partial^{<\underline i>}$ is the dual basis to $\underline \tau^{\{\underline i\}}$ and multiplication on differentials is given by
$$
\underline \partial^{<\underline k>} \underline \partial^{<\underline l>} = \left<{\underline k+ \underline l \atop \underline k}\right> \partial^{< \underline k+ \underline l>} \quad \mathrm{with} \quad \left<{\underline k+ \underline l \atop \underline k}\right > = \frac {\left({\underline k+ \underline l \atop \underline k}\right)} {\left\{{\underline k+ \underline l \atop \underline k}\right\}}.
$$
We also have
$$
\underline \partial^{<\underline k>} f = \sum_{\underline i \leq \underline k} \left\{ {\underline k \atop \underline i}\right \} \underline \partial^{<\underline i>}(f) \underline \partial^{<\underline k - \underline i>}.
$$
In this last formula, we implicitly make $\mathcal D^{(m)}_{X}$ act on $\mathcal O_{X}$.
This is formally obtained as follows: a differential operator of order $n$ is nothing but a linear map $P : \mathcal P^{n}_{Xm} \to \mathcal O_{X}$ and we compose it on the left with the map induced by the second projection $p_{2}^* : \mathcal O_{X} \to \mathcal P^{n}_{Xm}$.

For example, if we work locally, then $\underline t^{\underline h}$ is sent by $p_{2}^*$ to
$$
(\underline t + \underline \tau)^{\underline h} = \sum_{\underline k} {\underline h \choose \underline k} t^{\underline h-\underline k} \tau^{\underline k} = \sum_{ \underline k} q_{\underline  k}! {\underline  h \choose \underline k} t^{\underline h-\underline k} \tau^{\{\underline k\}}
$$
and therefore,
$$
\underline \partial^{< \underline k>}(\underline t^{\underline h}) = q_{\underline  k}!{\underline  h \choose \underline k} t^{\underline h-\underline k}.
$$
Finally, note that $\mathcal D^{(0)}_{X}$ is locally generated by $\partial_{1}, \ldots, \partial_{r}$ and that $\mathcal D^{(\infty)}_{X}$ is Grothendieck's ring of differential operators.
In general, when $k < p^{m+1}$, it is convenient to define $\partial^{[k]} = \partial^{<k>}/q_{k}!$ and note that $\mathcal D^{(m)}_{X}$ is locally generated by the $\partial_{i}^{[p^l]} = \partial^{<p^l>}$ for $l \leq m$.
In particular, we see that the diamond brackets notation should not appear very often in practice.

\section{The $p^m$-curvature map} \label{pem}

We assume from now on that $m \neq \infty$.

If $X$ is a scheme of characteristic $p$, we will denote by $F : X \to X$ the $m+1$-st iterate of its Frobenius endomorphism (given by the identity on $X$ and the map $f \mapsto f^{p^{m+1}}$ on functions).

\begin{lem}\label{frobgen}
Let $X \hookrightarrow Y$ be an immersion defined by an ideal $\mathcal I$. Then, the map
$$
\xymatrix@R=0cm{\mathcal I \ar[r] & \mathcal P_{Xm}(Y) \\ \varphi \ar@{|->}[r] & \varphi^{\{p^{m+1}\}},}
$$
composed with the projection
$$
\mathcal P_{Xm}(Y) \mapsto \mathcal P_{Xm}(Y)/\mathcal I \mathcal P_{Xm}(Y),
$$
is an $F^*$-linear map that is zero on ${\mathcal I}^2$.
\end{lem}

\textbf{Proof :}
If $\varphi, \psi \in \mathcal I$, we have
$$
(\varphi + \psi)^{\{p^{m+1}\}} = \varphi^{\{p^{m+1}\}} + \psi^{\{p^{m+1}\}} + \sum_{i + j = p^{m+1} \atop i,j >0} \left < {p^{m+1}} \atop i \right > \varphi^{\{i\}}\psi^{\{j\}}.
$$
When $0 < i, j < p^{m+1}$, we have $q_{i}, q_{j} < p$ and $q_{i}!$ and $q_{j}!$ are therefore invertible.
It follows that the last part in the sum falls inside $\mathcal I^{p^{m+1}}$.
In particular, it is zero modulo $\mathcal I \mathcal P_{Xm}(Y)$ and it follows that the composite map is additive.

Also, clearly, if $f \in \mathcal O_{X}$ and if $\varphi \in \mathcal I$, then $f\varphi$ is sent to
$$
(f\varphi)^{\{p^{m+1}\}} = f^{p^{m+1}} \varphi^{\{p^{m+1}\}} = F^{*} (f) \varphi^{\{p^{m+1}\}}.
$$
And we see that the map is $F^*$-linear.
Finally, if $\varphi, \psi \in \mathcal I$, then $\varphi\psi$ is sent to
$$
(\varphi\psi)^{\{p^{m+1}\}} = \varphi^{p^{m+1}} \psi^{\{p^{m+1}\}} \in \mathcal I \mathcal P_{Xm}(Y).
\quad \Box$$

For the rest of this section, we fix a base scheme $S$ of characteristic $p$ and we assume that $X$ is an $S$-scheme. We consider the usual commutative diagram with cartesian square (recall that here $F$ denotes the $m+1$-st iteration of Frobenius)
$$
\xymatrix{
X 
\ar[r]_{F_{X}} \ar@/_/[dr] \ar@/^1pc/[rr]^{F}
& X'  \ar[r] \ar[d] & X \ar[d]
\\
& S \ar[r]^F & S.
}
$$

If we apply Lemma \ref{frobgen} to the case of the diagonal embedding of $X$ in $X \times_{S} X$, we obtain, after linearizing and since $F_{X}^*\Omega^1_{X'} = F^*\Omega^1_{X}$, an $\mathcal O_{X}$-linear map
$$
\xymatrix@R=0cm{F_{X}^*\Omega^1_{X'} \ar[r] & \mathcal \mathcal P_{Xm}/\mathcal I \mathcal P_{Xm},}
$$
that we will call \emph{divided Frobenius}.

We may now prove the level $m$ version of Mochizuki's theorem (\cite{OgusVologodsky07}, Proposition 1.7):

\begin{prop} If $X$ is a smooth $S$-scheme, the divided Frobenius extends uniquely to an isomorphism of $\mathcal O_{X}$-modules
$$
F_{X}^*\Omega^1_{X'/S} \simeq \mathcal I_{Xm}\mathcal P^{p^{m+1}}_{Xm}/\mathcal I \mathcal P^{p^{m+1}}_{Xm}.
$$
\end{prop}

\textbf{Proof :}
From the discussion above, it is clear that we have such a map.
In order to show that this is an isomorphism, we may assume that there are local coordinates $t_{1}, \ldots, t_{r}$ on $X$, pull them back as $t'_{1}, \ldots, t'_{r}$ on $X'$ and also set as usual $\tau_{i} := 1 \otimes t_{i} - t_{i} \otimes 1$.
Then, our map is simply
$$
\xymatrix@R=0cm{\bigoplus_{i=1}^r \mathcal O_{X} \mathrm dt'_{i} \ar[r]^-\simeq & \bigoplus_{i=1}^r \mathcal O_{X} \tau_{i}^{\{p^{m+1}\}} \\ \mathrm dt'_{i} \ar@{|->}[r] & \tau_{i}^{\{p^{m+1}\}}. & \Box}
$$

Actually, we can do a little better.

\begin{prop} \label{pdfrob} If $X$ is a smooth $S$-scheme, then $\mathcal I_{Xm}^{\{p^{m+1}\}} \cap \mathcal I \mathcal P_{Xm}$ is stable under usual divided powers.
Moreover, the divided Frobenius extends uniquely to an isomorphism of divided power $\mathcal O_{X}$-algebras
$$
F_{X}^* \Gamma_{\bullet}\Omega^1_{X'/S} \simeq \mathcal P_{Xm}/\mathcal I \mathcal P_{Xm}.
$$
\end{prop}

\textbf{Proof :}
The first question is local and we assume for the moment that it is solved.
Then, by definition of the divided power algebra on a module, the above map
$$
F_{X}^*\Omega^1_{X'} \to \mathcal P_{Xm}/\mathcal I \mathcal P_{Xm}
$$
extends uniquely to a morphism of divided power algebras
$$
F_{X}^{*} \Gamma_{\bullet}\Omega^1_{X'} = \Gamma_{\bullet}F_{X}^{*}\Omega^1_{X'} \to \mathcal P_{Xm}/\mathcal I \mathcal P_{Xm}.
$$
Showing that it is an isomorphism is local again.

Thus, we assume that there are local coordinates $t_{1}, \ldots, t_{r}$ on $X$, we pull them back as $t'_{1}, \ldots, t'_{r}$ on $X'$ and we also set as usual $\tau_{i} := 1 \otimes t_{i} - t_{i} \otimes 1$.

We have
$$
\mathcal P_{Xm}/\mathcal I \mathcal P_{Xm} = \mathcal O_{X}<\tau_{1}, \ldots, \tau_{r} >^{(m)}/(\tau_{1}, \ldots, \tau_{r})
$$
which is therefore a free $\mathcal O_{X}$-module with basis $\underline \tau^{\{\underline kp^{m+1}\}}$.
The first assertion easily follows.

Moreover, our map is
$$
\xymatrix@R=0cm{\mathcal O_{X}<\mathrm dt'_{1}, \ldots, \mathrm d t'_{r} >^{(0)} \ar[r] & \mathcal O_{X}<\tau_{1}, \ldots, \tau_{r} >^{(m)}/(\tau_{1}, \ldots, \tau_{r})  \\ \mathrm dt'_{i} \ar@{|->}[r] & \tau_{i}^{\{p^{m+1}\}}.}
$$
And the left hand side is the free $\mathcal O_{X}$-module with basis $\underline {\mathrm dt'}^{[\underline k]}$.
Our assertion is therefore a consequence of the first part of Lemma \ref{compd} below.
$\Box$

\begin{lem} \label{compd}
In an ideal with partial divided powers of level $m$, we always have
\begin{enumerate}
\item For any $k \in \mathbf N$,
$$
(f^{\{p^{m+1}\}})^{[k]} = \frac {(kp)!}{(p!)^kk!} f^{\{kp^{m+1}\}}
$$
and $\frac {(kp)!}{(p!)^kk!} \in 1 + p\mathbf Z$.
\item If $t = qp^{m}+r$ with $q < p$ and $r < p^{m}$, then
$$
f^{\{kp^{m+1}\}}f^{\{t\}} = {kp+q \choose q} f^{\{kp^{m+1}+t\}}
$$
and ${kp+q \choose q} \in 1 + p\mathbf Z$.
\end{enumerate}
\end{lem}

\textbf{Proof :}
The first assertion comes from the case $m=0$ applied to $f^{p^m}$.
And we may consider the formula
$$
(f^{[p]})^{[k]} = \frac {(kp)!}{(p!)^kk!} f^{[kp]}
$$
as standard.
Moreover, there exists a product formula for the factor:
$$
\frac {(kp)!}{(p!)^kk!} = \prod_{j=1}^{k-1} {jp + p -1 \choose p-1}
$$
and it is therefore sufficient to prove that each factor in this product falls into  $1 + p\mathbf Z$.
We already know that they belong to $\mathbf Z$ and we have the product formula in $\mathbf Z_{(p)}$:
$$
{jp + p -1 \choose p-1} = \prod_{i=1}^{p-1}(1 + \frac ji p).
$$

The second assertion is even easier and comes from
$$
f^{\{u\}}f^{\{t\}} = \left\{{u \atop t}\right\} f^{\{u+t\}} \quad \mathrm{and} \quad  \left\{{kp^{m+1}+t \atop t}\right\} = {kp+q \choose q} = \prod_{i=1}^{q}(1 + \frac ki p).
\quad \Box$$

\begin{dfn}
If $X$ is a smooth $S$-scheme, the \emph{$p^m$-curvature} map is the morphism
$$
F_{X}^*S^{\bullet} \mathcal T_{X'} \to \mathcal D_{X}^{(m)}
$$
obtained by duality from the composite
$$
\mathcal P_{Xm} \to \mathcal P_{Xm}/\mathcal I \mathcal P_{Xm} \simeq F_{X}^{*}\Gamma_{\bullet}\Omega^1_{X'}
$$
\end{dfn}

We have to be a little careful here : first of all, we consider the induced morphisms
$$
\mathcal P^{[k]}_{Xm} \to F_{X}^{*}\Gamma_{\leq k}\Omega^1_{X'},
$$
(with usual divided powers on $\mathcal P_{Xm}$), and then we dualize to get
$$
F_{X}^{*}S^{\leq k} \mathcal T_{X'} \to \mathcal D_{X,k}^{(m)}.
$$
and take the direct limit on both sides.

Alternatively, we may also call \emph{$p^m$-curvature} the adjoint map
$$
S^{\bullet} \mathcal T_{X'} \to F_{X*}\mathcal D_{X}^{(m)}.
$$

We will denote by $\mathcal Z_{X}^{(m)}$ the center of $\mathcal D_{X}^{(m)}$ and by $\mathcal Z \mathcal O_{X}^{(m)}$ the centralizer of $\mathcal O_{X}$ in $\mathcal D_{X}^{(m)}$.

\begin{prop}
If $X$ is a smooth scheme over $S$, then the $p^m$-curvature map induces an isomorphism of $\mathcal O_{X}$-algebras $F_{X}^*S^{\bullet} \mathcal T_{X'} \simeq \mathcal Z\mathcal O_{X}^{(m)}$ and an isomorphism of $\mathcal O_{X'}$-algebras $S^{\bullet} \mathcal T_{X'} \simeq F_{X*}\mathcal Z_{X}^{(m)}$.
\end{prop}

Note that it will formally follow from the definition of the multiplication in $\mathcal D_{X}^{(m)}$ and Proposition \ref{dualmul} that the $p^m$-curvature map is a morphism of algebras.
More precisely, this map is obtained by duality from a morphism of coalgebras.
However, we need a local description in order to prove the rest of the proposition.

\textbf{Proof :}
Both questions are local and we may therefore use local coordinates $t_{1}, \ldots, t_{r}$, pull them back to $t'_{1}, \ldots, t'_{r}$ on $X'$ and denote by $\xi'_{1}, \ldots, \xi'_{r}$ the corresponding basis of $\mathcal T_{X'}$.
By construction, the $p^{m}$-curvature map is then given by
$$
{\xi'}_{i}^k \mapsto  \partial_{i}^{<kp^{m+1}>}.
$$
It follows from Lemma \ref{compd} (and duality) that
$$
\partial_{i}^{<kp^{m+1}>} = (\partial_{i}^{<p^{m+1}>})^{k}
$$
and this shows that we do have a morphism of rings, which is clearly injective because we have free modules on both sides.
The image is the $\mathcal O_{X'}$-subalgebra generated by the $\partial_{i}^{<p^{m+1}>}$ and this is exactly the center as Berthelot showed in Proposition 2.2.6 of \cite{Berthelot96}.
$\Box$

The following theorem is due to Masaharu Kaneda (\cite{Kaneda04}, section 2.3, see also \cite{BezrukavnikovMirkovicRumynin08} in the case $m = 0$) when $S$ is the spectrum of an algebraically closed field.

\begin{thm}
Let $X$ be a smooth scheme over a scheme $S$ of positive characteristic $p$ and $F_{X} : X \to X'$ the $m+1$-st iterate of the relative Frobenius.
Let $\mathcal D^{(m)}_{X}$ be the ring of differential operators of level $m$ on $X/S$ and $\mathcal Z\mathcal O^{(m)}_{X}$ the centralizer of $\mathcal O_{X}$.
Then, there is an isomorphism of $\mathcal Z\mathcal O^{(m)}_{X}$-algebras
$$
\xymatrix@R=0cm{F_{X}^*F_{X*} \mathcal D^{(m)}_{X} \ar[r] &\mathcal End_{\mathcal Z\mathcal O^{(m)}_{X}}(\mathcal D^{(m)}_{X}) \\ f \otimes Q \ar@{|->}[r]  & (P \mapsto fPQ).}
$$
\end{thm}

\textbf{Proof :}
The question is local and one easily sees that $\mathcal D_{X/S}^{(m)}$ is free as $\mathcal Z\mathcal O^{(m)}_{X}$-modules on the generators $\underline \partial^{<\underline k>}$ with $\underline k < \underline {p^{m+1}}$.
More precisely, this follows again from Lemma \ref{compd} that gives us, by duality,
$$
\partial_{i}^{<kp^{m+1}+t>} = (\partial_{i}^{<p^{m+1}>})^{k} \partial_{i}^{<t>} 
$$
when $t < p^{m+1}$.
It is then sufficient to compare basis on both sides (see Kaneda's proof for the details).
$\Box$

For example, when $m =0$, in the simplest case of an affine curve $X = \mathrm{Spec} A$ with coordinate $t$ and corresponding derivation $\partial$, the first powers $1, \partial, \ldots, \partial^{p-1}$ form a basis of the ring of differential operators $\mathcal D$ and the map of the theorem sends $\partial^k$ to
$$
\left[\begin{array}{cc}0 & \partial^pI_{k} \\I_{p-k} & 0\end{array}\right] \in M_{p\times p}(A[\partial^p])
$$
for $k = 0, \ldots, p-1$, and $\partial^p$ to $\partial^pI_{p}$.

This theorem is usually stated as \emph{proving the Azumaya nature of $F_{X*}\mathcal D_{X}^{(m)}$}.
More precisely, we can see $F_{X*}\mathcal D_{X}^{(m)}$ as a sheaf of algebras on
$$
\check{\mathbf T}_{X'} = \mathrm{Spec} S^{\bullet} \mathcal T_{X'} \simeq \mathrm{Spec} F_{X*}\mathcal Z_{X}^{(m)}
$$
and the above theorem provides a trivialization of $F_{X*}\mathcal D_{X}^{(m)}$ along the ``Frobenius''
$$
F_{X} : X \otimes_{X'} \check{\mathbf T}_{X'} \to \check{\mathbf T}_{X'}.
$$

\begin{prop} Let $X$ be a smooth $S$-scheme.
If we denote by $\mathcal K_{X}^{(m)}$ the two-sided ideal of $\mathcal D_{X}^{(m)}$ generated by the image of $ \mathcal T_{X'}$ under the $p^m$-curvature map, there is an exact sequence
$$
0 \to \mathcal K_{X}^{(m)} \to \mathcal D_{X}^{(m)} \to \mathcal End_{\mathcal O_{X'}}(\mathcal O_{X})\to 0.
$$
\end{prop}

\textbf{Proof :} This follows from \cite{Berthelot96}, Proposition 2.2.7.
$\Box$

We will denote by  $\widehat {\mathcal D}_{X}^{(m)}$ the completion of $\mathcal D_{X}^{(m)}$ along the two-sided ideal $\mathcal K_{X}^{(m)}$.
We will also denote by $\widehat {\mathcal Z}_{X}^{(m)}$ the completion of $\mathcal Z_{X}^{(m)}$ along $\mathcal Z_{X}^{(m)} \cap \mathcal K_{X}^{(m)}$ and by $\widehat{\mathcal Z}\mathcal O^{(m)}_{X}$ the completion of $\mathcal Z\mathcal O^{(m)}_{X}$ along $\mathcal Z\mathcal O^{(m)}_{X} \cap \mathcal K_{X}^{(m)}$.
Note that the $p^m$-curvature map gives isomorphisms
$$
F^*_{X}\widehat {S^{\bullet} \mathcal T_{X'}} \simeq \widehat {\mathcal Z}\mathcal O_{X}^{(m)} \quad \mathrm{and} \quad  \widehat {S^{\bullet} \mathcal T_{X'}} \simeq F_{X*}\widehat {\mathcal Z}_{X}^{(m)}.
$$

\begin{prop} If $X$ is a smooth $S$-scheme, we have natural isomorphisms
$$
\widehat {\mathcal Z}_{X}^{(m)} \otimes_{\mathcal Z_{X}^{(m)}} \mathcal D_{X}^{(m)} \simeq \widehat {\mathcal D}_{X}^{(m)} \simeq \mathcal Hom_{\mathcal O_{X}}(\mathcal P_{Xm}, \mathcal O_{X}).
$$
\end{prop}

\textbf{Proof :} The existence of the first map is clear and it formally follows from the definitions that it is an isomorphism.

Now, note that the canonical projections $\mathcal P_{Xm} \to \mathcal P^n_{Xm}$ induce a compatible family of maps
$$
\mathcal D_{Xn}^{(m)} = \mathcal Hom_{\mathcal O_{X}}(\mathcal P^n_{Xm}, \mathcal O_{X}) \to \mathcal Hom_{\mathcal O_{X}}(\mathcal P_{Xm}, \mathcal O_{X})
$$
from which we derive a morphism $\mathcal D_{X}^{(m)} \to \mathcal Hom_{\mathcal O_{X}}(\mathcal P_{Xm}, \mathcal O_{X})$.
On the other hand, using the isomorphism of Proposition \ref{pdfrob} and the $p^m$-curvature, the projection $\mathcal P_{Xm} \to \mathcal P_{Xm}/\mathcal I \mathcal P_{Xm}$ dualizes to
$$
\widehat {\mathcal Z}_{X}^{(m)} \hookrightarrow  \mathcal Hom_{\mathcal O_{X}}(\mathcal P_{Xm}, \mathcal O_{X}).
$$
And by construction, these two maps are compatible on $\mathcal Z_{X}^{(m)}$ and induce a map
$$
\widehat {\mathcal Z}_{X}^{(m)} \otimes_{\mathcal Z_{X}^{(m)}} \mathcal D_{X}^{(m)} \to \mathcal Hom_{\mathcal O_{X}}(\mathcal P_{Xm}, \mathcal O_{X}).
$$
It is now a local question to check that this is an isomorphism.
$\Box$

\section{Lifting the $p^m$-curvature} \label{lift}

We will prove here the Azumaya nature of the ring of differential operators of higher level.
In order to make it easier to read, we will not always mention direct images under Frobenius.
This is not very serious because Frobenius maps are homeomorphisms and playing with direct image only impacts the linearity of the maps (and this should be clear from the context).

\begin{dfn}
If $X$ is a scheme of characteristic $p$, a \emph{lifting} $\widetilde X$ of $X$ modulo $p^2$ is a flat $\mathbf Z/p^2\mathbf Z$-scheme $\widetilde X$ such that $X = \widetilde X \times_{\mathbf Z/p^2\mathbf Z} \mathbf F_{p}$.
A \emph{lifting of a morphism} $f : Y \to X$ of schemes of characteristic $p$ is a morphism $\widetilde f : \widetilde X \to \widetilde Y$ between liftings such that $f = \widetilde f \times_{\mathbf Z/p^2\mathbf Z} \mathbf F_{p}$.
\end{dfn}

We will use the well known elementary result:

\begin{lem}
If $M$ is a $\mathbf Z/p^2\mathbf Z$-module, multiplication by $p!$ induces a surjective map
$$
p! : M/pM \to pM
$$
which is bijective if $M$ is flat.
\end{lem}

\textbf{Proof :} Exercise. $\Box$

Note that $p! = -p \mod p^2$ and this explains why minus signs will appear in the formulas below.
Actually, we will need more fancy estimates :

\begin{lem}\label{lucas}
We have for any $m > 0$,
$$
{{p^{m+1}} \choose i} = \left\{  \begin{array}{cl} 1 & \mathrm{if} \quad i=0 \quad \mathrm{or} \quad i= p^{m+1} \\ (-1)^k p! & \mathrm{if} \quad i = kp^m \\ 0 & \mathrm{otherwise} \end{array} \right. \mod p^2.
$$
\end{lem}

\textbf{Proof :} Standard results on valuations of factorials show that
$$
v_{p}({p^{m+1} \choose i}) =  \left\{\begin{array} {lll} 0 & \mathrm{if} & i = 0 \quad \mathrm{or} \quad i = p^{m+1} \\ 1 & \mathrm{if} & i = kp^m \quad \mathrm{with} \quad 0 < k < p \\ > 1 & \mathrm{otherwise}  \end{array} \right.
$$
and we are therefore reduced to showing that
$$
{{p^{m+1}} \choose  kp^m} =  (-1)^k p!  \mod p^2,
$$
or what is slightly easier, that
$$
{{p^{m+1} - 1} \choose  kp^m - 1} =  (-1)^{k+1} \mod p.
$$
First of all, we can use Lucas congruences that give
$$
{{p^{m+1} - 1} \choose  kp^m - 1} = {{p - 1} \choose  k - 1} \mod p,
$$
and then the binomial property
$$
 {{p - 1} \choose  k - 1} + {{p - 1} \choose  k - 2} =  {p \choose  k - 2} = 0  \mod p
$$
in order to reduce to the case $k = 1$.
$\Box$

Up to the end of the section, we let $S$ be a scheme of characteristic $p$ and denote by $\widetilde S$ a lifting of $S$ (modulo $p^2$ as defined above).

\begin{dfn}
If $X$ is an $S$-scheme, a strong lifting $\widetilde F : \widetilde X \to \widetilde X'$ of the $m+1$-st iteration of Frobenius of $X$ is a morphism that satisfies
$$
f' = 1 \otimes f \mod p \quad \Rightarrow \quad  \widetilde F^*(f') = f^{p^{m+1}} + pg^{p^{m}} \quad \mathrm{with} \quad  \mathrm g \in \mathcal O_{\widetilde X}.
$$
\end{dfn}

When $m = 0$, this is nothing but a usual lifting but the condition is stronger in general.
For example, the map $t \mapsto t^4+2t$ is not a strong lifting of the Frobenius on the affine line when $m =1$ and $p = 2$.
However, the condition is usually satisfied in practice, especially when the lifting comes from a lifting of the absolute Frobenius as the next lemma shows.

\begin{lem}
If $\widetilde F : \widetilde X \to \widetilde X$ is a lifting of the \emph{true absolute} Frobenius of $X$, then for $f \in \mathcal O_{\widetilde X}$, we have
$$
\widetilde F^{m+1*}(f) = f^{p^{m+1}} + pg^{p^{m}}
$$
with $g \in \mathcal O_{\widetilde X}$.
\end{lem}

\textbf{Proof :}
By definition, we  can write
$$
\widetilde F^*(f) = f^{p} + pg
$$
with $g \in \mathcal O_{\widetilde X}$ and we prove by induction on $m$ that
$$
\widetilde F^{m+1*}(f) = f^{p^{m+1}} + pg^{p^{m}}.
$$
If we apply the ring homomorphism $\widetilde F^*$ on both sides of this equality, we get
$$
\widetilde F^{m+2*}(f) = \widetilde F^{*}(f)^{p^{m+1}} + p\widetilde F^{*}(g)^{p^{m}}
$$
$$
= (f^p + pg)^{p^{m+1}} + p(g^p + ph)^{p^{m}} = f^{p^{m+2}} + pg^{p^{m+1}}.
\quad \Box$$

Now, we fix a smooth $S$-scheme $X$ and let $\widetilde F : \widetilde X \to \widetilde X'$ be a strong lifting of the $m+1$-st iteration of the relative Frobenius of $X$.
We will denote by $\widetilde X \times \widetilde X$  (resp. $\widetilde X' \times \widetilde X'$) the fibered product over $\widetilde S$ and by $\widetilde {\mathcal I}$ (resp. $\widetilde {\mathcal I}'$) the ideal of $\widetilde X$ in $\widetilde X \times \widetilde X$  (resp. $\widetilde X'$ in $\widetilde X' \times \widetilde X'$).

\begin{lem} \label{compute}
Assume that $\widetilde F^*(f') = f^{p^{m+1}} + pg^{p^{m}}$ with $g \in \mathcal O_{\widetilde X}$.
Let $\varphi = 1 \otimes f - f \otimes 1$, $\varphi' = 1 \otimes f' - f' \otimes 1$ and $\psi = 1 \otimes g - g \otimes 1$.
Then, the composite map
$$
\widetilde F^* : \widetilde {\mathcal I'} \hookrightarrow  \mathcal O_{\widetilde X' \times \widetilde X'} \stackrel {\widetilde F^* \times \widetilde F^* }\longrightarrow \mathcal O_{\widetilde X \times \widetilde X} \to \mathcal P_{\widetilde X,m}.
$$
sends $\varphi' $ to
$$
p!\left(\varphi^{\{p^{m+1}\}} +  \sum_{k=1}^{p-1} (-1)^k f^{(p-k)p^{m}} \varphi^{kp^m} - \psi^{p^m} \right).
$$
\end{lem}

\textbf{Proof :}
We have
$$
\widetilde F^*(\varphi') = 1 \otimes f^{p^{m+1}} - f^{p^{m+1}} \otimes 1 + 1 \otimes pg^{p^m} - pg^{p^m} \otimes 1
$$
(recall that $p^2 = 0$ on $\widetilde S$)
$$
= \varphi^{p^{m+1}} + \sum_{i=1}^{p^{m+1}-1}{p^{m+1} \choose i} f^{p^{m+1}-i}\varphi^i + p\psi^{p^m}.
$$
We finish with Lemma \ref{lucas}.
$\Box$

\begin{prop} \label{liffrob}
There is a well defined map
$$
\frac 1{p!} \widetilde F^* :  \widetilde {\mathcal I'} \to p\mathcal P_{\widetilde X,m} \simeq \mathcal P_{\widetilde X,m}/p\mathcal P_{\widetilde X,m} \simeq \mathcal P_{X,m}
$$
that factors through $\Omega^1_{X'}$ and takes values inside $\mathcal J_{Xm}$.
Moreover, the induced morphism
$$
\frac 1{p!} \widetilde F^* :  \Omega^1_{X'} \to \mathcal P_{X,m}
$$
is a lifting of divided Frobenius
$$
\Omega^1_{X'} \to  \mathcal P_{X,m}/\mathcal I\mathcal P_{X,m}.
$$
\end{prop}

\textbf{Proof :}
It follows from lemma \ref{compute} that the map is well defined: more precisely, we need to check that $\widetilde F^*$ sends $\widetilde {\mathcal I'}$ inside $p\mathcal P_{\widetilde X,m}$.
By linearity, it is sufficient to consider the action on sections $\varphi'$ as in the lemma.

Now, since $\widetilde F^*$ is a morphism of rings that sends $\widetilde {\mathcal I'}$ to zero modulo $p$, it is clear that $\frac 1{p!} \widetilde F^*$ will send $\widetilde {\mathcal I'}^2$ to $0$.
Thus, it factors through $\Omega^1_{\widetilde X'}$.
Actually, since the target is killed by $p$, it even factors through $\Omega^1_{X'}$.
And it falls inside $\mathcal J_{Xm}$ thanks to the first part.

Finally, the last assertion follows again from the explicit description of the map.
$\Box$

\textbf{Warning} : The quotient map $\mathcal P_{X,m} \to \mathcal P_{X,m}/\mathcal I\mathcal P_{X,m}$ is \emph{not} compatible with the divided power structures : $\tau_{i}^{p^m}$ is sent to $0$ but $(\tau_{i}^{p^m})^{[p]} = \tau_{i}^{\{p^{m+1}\}}$ is not.

\begin{prop}
The divided Frobenius $\frac 1{p!} \widetilde F^*$ extends canonically to a morphism
$$
F_{X}^*\Gamma_{\bullet}\Omega^1_{X'} \to  \mathcal P_{X,m}.
$$
By duality, we get a morphism of $\mathcal O_{X}$-modules
$$
\Phi_{X} : \widehat {\mathcal D}_{X}^{(m)} \longrightarrow \widehat{\mathcal Z}\mathcal O^{(m)}_{X} \hookrightarrow\widehat {\mathcal D}_{X}^{(m)}.
$$
\end{prop}

\textbf{Proof :} We saw in Proposition \ref{liffrob} that the morphism $\frac 1{p!} \widetilde F^*$ takes values into $\mathcal J_{X,m}$ and therefore extends to a morphism of divided power algebras
$$
\Gamma_{\bullet}\Omega^1_{X'} \to  \mathcal P_{X,m}
$$
that we can linearize.
Moreover, we saw in Proposition \ref{compd} that
$$
\widehat {\mathcal D}_{X}^{(m)} \simeq \mathcal Hom_{\mathcal O_{X}}(\mathcal P_{Xm}, \mathcal O_{X})
$$
and we also have
$$
\widehat{\mathcal Z}\mathcal O^{(m)}_{X} \simeq F_{X}^*\widehat{S^{\bullet}\mathcal T_{X'}} \simeq F_{X}^*\mathcal Hom_{\mathcal O_{X'}}(\Gamma_{\bullet}\Omega^1_{X'}, \mathcal O_{X'}) \simeq \mathcal Hom_{\mathcal O_{X}}(F_{X}^*\Gamma_{\bullet}\Omega^1_{X'}, \mathcal O_{X}).
\quad \Box$$

\begin{dfn}
The morphism $\Phi_{X}$ is called the \emph{Frobenius} of $\widehat {\mathcal D}_{X}^{(m)}$.
\end{dfn}

We will simply write $\Phi$ when $X$ is understood from the context but we might also use $\Phi_{X}^{(m)}$ to indicate the level.
Note that $\Phi$ actually depends on the choice of the strong lifting $\widetilde F$ of $F_X$.

\begin{prop}\label{locdesc}
If we are given local coordinates $t_{1}, \ldots, t_{r}$, then
$$
\Phi(\underline \partial^{\langle \underline n \rangle}) =
\left\{\begin{array}{l l}
1 & \mathrm{if} \quad |\underline n| = \underline 0
\\
0 & \mathrm{if} \quad 0 < |\underline n| < p^m
\\
 \frac 1{p!} \sum_{j=1}^r \partial_{i}^{[p^m]}(\widetilde F^*(\tilde {t'_{j}})) \partial_{j}^{\langle p^{m+1}\rangle } & \mathrm{if} \quad \underline n = p^m\mathbf 1_{i}.
\end{array}\right.
$$
Actually, if we have $\widetilde F^*(\widetilde {t'_{j}}) = \widetilde {t_{j}}^{p^{m+1}} + p\widetilde {g_{j}}
^{p^m}$, the third expression can be rewritten
$$
- t_{i}^{(p-1)p^{m}}\partial_{i}^{\langle p^{m+1}\rangle } - \sum_{j=1}^r \partial_{i}(g_{j})^{p^m}  \partial_{j}^{\langle p^{m+1}\rangle } .
$$
\end{prop}

\textbf{Proof :} 
The point consists in writing $\Phi(\underline \partial^{\langle \underline n \rangle })$ in the topological $\mathcal O_{X}$-basis $\underline \partial^{\langle \underline k p^{m+1}\rangle }$ of $\widehat{\mathcal Z}\mathcal O^{(m)}_{X}$ when $|\underline n| \leq p^m$.
By duality, the coefficient of $\underline \partial^{\langle \underline k p^{m+1}\rangle }$ in $\Phi(\underline \partial^{\langle \underline n\rangle })$ is identical to the coefficient of $\underline \tau^{\{\underline n\}}$ in the image of $1 \otimes (\underline{\mathrm{d}t})^{[\underline k]}$ under the morphism
$$
\mathcal O_{X} \otimes_{\mathcal O_{X'}}\Gamma_{\bullet}\Omega^1_{X'} \to  \mathcal P_{X,m}.
$$
Recall that this is exactly
$$
\left(\frac{ \widetilde F^*}{p!}(\widetilde {\underline \tau'})\right)^{[\underline k]}.
$$
Since we consider only the case $|\underline n| \leq p^m$ we may work modulo $\mathcal I^{\{p^m+1\}}$ on the right.
If we write $\widetilde F^*(\widetilde {t'_{j}}) = \widetilde {t_{j}}^{p^{m+1}} + p\widetilde {g_{j}}
^{p^m}$, we obtain
$$
\frac 1{p!}\widetilde F^*(\widetilde \tau'_{j}) = \frac 1{p!}\ \sum_{\underline s \neq 0} \underline \partial^{\langle \underline s \rangle}(\widetilde F^*(\widetilde t'_{j})) \tau^{\{\underline s\}}
$$
$$
=  - t_{j}^{(p-1)p^{m}} \tau_{j}^{p^m} - \sum_{i=1}^r \partial_{i}(g_{j})^{p^m} \tau_{i}^{p^m} \mod {\mathcal I}^{\{p^m+1\}}.
$$
Thus, we see that the only contributions will come from the case $|\underline k| \leq 1$ and that, when $\underline k \neq \underline 0$, the coefficient of $\underline \tau^{\{\underline n\}}$ is zero unless $\underline n = p^m\mathbf 1_{i}$.
Then, there are two cases, first $i\neq j$ in which case, only
$$
 \frac 1{p!}\partial_{i}^{[p^m]}(\widetilde F^*(\tilde {t'_{j}})) = -\partial_{i}(g_{j})^{p^m}
$$ is left and the case $i=j$ where we obtain
$$
 \frac 1{p!} \partial_{i}^{[p^m]}(\widetilde F^*(\tilde {t'_{i}})) = - t_{i}^{(p-1)p^{m}} - \partial_{i}(g_{i})^{p^m}.
\quad \Box$$

For example, when $m = 0$, in the case of  the affine line with parameter $t$ and derivation $\partial$ , if we choose the usual lifting of Frobenius $t \mapsto t^p$, we obtain the simple formula $\Phi(\partial) = - t^{p-1}\partial^p$.

Formulas are a lot more complicated in general but they become surprisingly nice when we stick to the usual generators of the center.

\begin{prop}\label{phibar}
For all $i = 1, \cdots, r$, we have
$$
\Phi(\partial_{i}^{\langle p^{m+1} \rangle}) = \partial_{i}^{\langle p^{m+1} \rangle} + \Phi(\partial_{i}^{\langle p^{m} \rangle})^p.
$$
\end{prop}

\textbf{Proof :} 
As above, the coefficient of $\underline \partial^{\langle \underline k p^{m+1}\rangle }$ in $\Phi(\partial_{i}^{\langle p^{m+1} \rangle})$ is identical to the coefficient of $\tau_{i}^{\{p^{m+1}\}}$ in
$$
\frac{ \widetilde F^*}{p!}(\widetilde {\underline \tau'})^{[\underline k]} = \prod_{j=1}^r\left(\tau_{j}^{\{p^{m+1}\}} +  \sum_{l=1}^{p-1}(1)^l t_{j}^{(p-l)p^{m}} \tau_{j}^{lp^m} - \sum_{l=1}^r \partial_{l}(g_{j})^{p^m} \tau_{l}^{p^m} \right)^{[k_{j}]}
$$
if we write $\widetilde F^*(\widetilde {t'_{j}}) = \widetilde {t_{j}}^{p^{m+1}} + p\widetilde {g_{j}}^{p^m}$.

Thus, we see that the only contributions will come from the cases $\underline k =\mathbf 1_{i}$ that gives $\tau_{i}^{\{p^{m+1}\}}$ and $\underline k = p\mathbf1_{j}$ that will give
$$
\left(- \partial_{i}(g_{j})^{p^m} \tau_{i}^{p^m}\right)^{[p]} = \left(- \partial_{i}(g_{j})^{p^m}\right)^{p} \tau_{i}^{\{p^{m+1}\}}
$$
for all $j$ plus the special contribution
$$
\left(- t_{i}^{(p-1)p^{m}} \tau_{i}^{p^m}\right)^{[p]} = \left(- t_{i}^{(p-1)p^{m}}\right)^{p}  \tau_{i}^{\{p^{m+1}\}}
$$
of the case $j = i$.
In other words, we obtain
$$
\Phi(\partial_{i}^{\langle p^{m+1} \rangle}) =  \partial_{i}^{\langle p^{m+1} \rangle} + \left(- t_{i}^{(p-1)p^{m}}\right)^p\partial_{i}^{\langle p^{m+2}\rangle } + \sum_{j=1}^r \left(-\partial_{i}(g_{j})^{p^m}\right)^p  \partial_{j}^{\langle p^{m+2}\rangle } .
$$
Our assertion therefore follows from the formulas in Proposition \ref{locdesc} because, thanks to Lemma \ref{compd}, we have
$$
\left(\partial_{j}^{\langle p^{m+1}\rangle }\right)^p =  \partial_{j}^{\langle p^{m+2}\rangle }.
\quad \Box$$

This calculation shows in particular that $\Phi$ is \emph{not} a morphism of rings.
However, we will see later on that the map induced by $\Phi$ on the center is a morphism of rings so that the above formula fully describes this map.

We recall now the following general result on Frobenius and divided powers:

\begin{lem} \label{frobfac}
The canonical map $P_{Xm} \to X \times X$ factors through $X \times_{X'} X$ if $X$ is seen as an $X'$-scheme via $F_{X}$.
Actually, if we are given local coordinates $t_{1}, \ldots, t_{r}$ and we set $\tau_{i} := 1 \otimes t_{i} - t_{i} \otimes 1$, the corresponding map on sections is the canonical injection
$$
\mathcal O_{X}[\underline \tau]/(\underline \tau^{p^{m+1}}) \hookrightarrow \mathcal O_{X}\langle \underline \tau\rangle^{(m)}.
$$
\end{lem}

\textbf{Proof :}
In order to prove the first assertion, the point is to check that if $f \in \mathcal O_{X}$, then $1 \otimes f^{p^{m+1}} - f^{p^{m+1}} \otimes 1$ is sent to zero in $\mathcal P_{Xm}$.
But we have
$$
1 \otimes f^{p^{m+1}} - f^{p^{m+1}} \otimes 1 = p! (1 \otimes f - f \otimes 1)^{\{p^{m+1}\}} = 0.
$$
Concerning the second assertion, we just have to verify that
$$
\mathcal O_{X}[\underline \tau]/(\underline \tau^{p^{m+1}}) \simeq \mathcal O_{X \times_{X'} X}.
$$
Since Frobenius is cartesian on \'etale maps (\cite{SGA5}, XIV, 1, Proposition 2), we may assume that $X = \mathbf A^r_{S}$, in which case this is clear.
$\Box$

Our main theorem arrives now:

\begin{thm} \label{mainth}
Let $X$ be a smooth scheme over a scheme $S$ of positive characteristic $p$ and $\widetilde F$ a strong lifting of the $m+1$-st iterate of the relative Frobenius of $X$.
The divided Frobenius extends canonically to an isomorphism
$$
\mathcal O_{X \times_{X'} X} \otimes_{\mathcal O_{X'}}\Gamma_{\bullet}\Omega^1_{X'} \simeq  \mathcal P_{X,m}.
$$
By duality, we obtain an isomorphism of $\mathcal O_{X}$-algebras
$$
\widehat {\mathcal D}_{X}^{(m)} \simeq \mathcal End_{\widehat {\mathcal Z}_{X}^{(m)}}(\widehat{\mathcal Z}\mathcal O^{(m)}_{X}).
$$
\end{thm}

\textbf{Proof :}
Thanks to lemma \ref{frobfac}, we may extend by linearity the morphism of divided power algebras
$$
\Gamma_{\bullet}\Omega^1_{X'} \to  \mathcal P_{X,m}
$$
and obtain
$$
\mathcal O_{X \times_{X'} X} \otimes_{\mathcal O_{X'}}\Gamma_{\bullet}\Omega^1_{X'} \to  \mathcal P_{X,m}.
$$
We now show that it is an isomorphism.
This is a local question and we may therefore fix some coordinates $t_{1}, \ldots, t_{r}$, call $t'_{1}, \dots, t'_{r}$ the corresponding coordinates on $X'$ and as usual, set $\tau_{i} := 1 \otimes t_{i} - t_{i} \otimes 1$.
It follows from proposition \ref{liffrob} that $\mathrm {d}t'_{i}$ is sent to $\tau_{i}^{\{p^{m+1}\}} + \phi_{i}$ with
$$
\phi_{i} \in (\underline \tau^{(p^m)})\mathcal O_{X}[\underline \tau]/(\underline \tau^{p^{m+1}}).
$$
And we have to check that the divided power morphism of $\mathcal O_{X}[\underline \tau]/(\underline \tau^{p^{m+1}})$-algebras
$$
\xymatrix@R=0cm{\mathcal O_{X}[\underline \tau]/(\underline \tau^{p^{m+1}})\langle \underline {\mathrm {d}t'} \rangle^{(0)} \ar[r] &\mathcal O_{X}\langle \underline \tau \rangle^{(m)} \\ \mathrm{d}t'_{i} \ar@{|->}[r] &  \tau_{i}^{\{p^{m+1}\}} + \phi_{i}.}
$$
is bijective.
This is easy because we have free modules with explicit basis on both sides.

Now, we obtain our assertion by duality.
Using the fact that $X$ is finite flat over $X'$, so that $\mathcal O_{X}$ is locally free over $\mathcal O_{X'}$, we have the following sequence of isomorphisms
$$
\mathcal Hom_{\mathcal O_{X}}(\mathcal O_{X \times_{X'} X} \otimes_{\mathcal O_{X'}}\Gamma_{\bullet}\Omega^1_{X'}, \mathcal O_{X}) \simeq \mathcal Hom_{\mathcal O_{X'}}(\Gamma_{\bullet}\Omega^1_{X'}, \mathcal Hom_{\mathcal O_{X}}(\mathcal O_{X \times_{X'} X}, \mathcal O_{X}))
$$
$$
\simeq \mathcal Hom_{\mathcal O_{X'}}(\Gamma_{\bullet}\Omega^1_{X'}, \mathcal End_{\mathcal O_{X'}}(\mathcal O_{X})) \simeq \mathcal Hom_{\mathcal O_{X'}}(\Gamma_{\bullet}\Omega^1_{X'}, \mathcal O_{X'}) \otimes_{\mathcal O_{X'}}\mathcal End_{\mathcal O_{X'}}(\mathcal O_{X})
$$
$$
\simeq \widehat{S^{\bullet}\mathcal T_{X'}} \otimes_{\mathcal O_{X'}}\mathcal End_{\mathcal O_{X'}}(\mathcal O_{X}) \simeq \mathcal End_{\widehat{S^{\bullet}\mathcal T_{X'}} }(\mathcal O_{X} \otimes_{\mathcal O_{X'}} \widehat{S^{\bullet}\mathcal T_{X'}} ).
$$
and we know that $S^{\bullet}\mathcal T_{X'} \simeq F_{X*}\mathcal Z_{X}^{(m)}$.

It remains to show that this is a morphism of rings and we do that by proving that it comes by duality from a morphism of coalgebras.
Actually, both morphisms $$ \mathcal O_{X \times_{X'} X} \to  \mathcal P_{X,m} \quad \mathrm{and} \quad \Gamma_{\bullet}\Omega^1_{X'} \to  \mathcal P_{X,m} $$ are compatible with the coalgebra structures.
We can be more precise: for the first one, this is because the comultiplication is induced on both sides by the same formula
$$
f \otimes g \mapsto f \otimes 1 \otimes g
$$
and for the second one, it is an immediate consequence of the universal property of divided powers.
$\Box$

\textbf{Warning}: The isomorphism of the theorem is \emph{not} a morphism of $\widehat {\mathcal Z}_{X}^{(m)}$-algebras.
However, we have the following:

\begin{cor}
The morphism $\Phi$ induces an automorphism of the ring $\widehat {\mathcal Z}_{X}^{(m)}$ (that depends on the lifting of Frobenius).
\end{cor}

\textbf{Proof :} The isomorphism of rings
$$
\widehat {\mathcal D}_{X}^{(m)} \simeq \mathcal End_{\widehat {\mathcal Z}_{X}^{(m)}}(\widehat{\mathcal Z}\mathcal O^{(m)}_{X}).
$$
induces an isomorphism on the centers which is nothing but the map induced by $\Phi$. $\Box$

\begin{cor}
The sheaf $\widehat{\mathcal Z}\mathcal O^{(m)}_{X}$ is a (left) $\widehat {\mathcal D}_{X}^{(m)}$-module for the action
$$
P \bullet fQ = \Phi(Pf)Q, \quad P \in \widehat {\mathcal D}_{X}^{(m)}, f \in \mathcal O_{X},Q \in \widehat {\mathcal Z}_{X}^{(m)}
$$
(that again depends on the lifting of Frobenius).
\end{cor}

\textbf{Proof :} Using the fact that the action is $\widehat{\mathcal Z}\mathcal O^{(m)}_{X}$-linear by definition, we may assume that $Q = 1$.
And since the isomorphism 
$$
\widehat {\mathcal D}_{X}^{(m)} \simeq \mathcal End_{\widehat {\mathcal Z}_{X}^{(m)}}(\widehat{\mathcal Z}\mathcal O^{(m)}_{X}).
$$
is $\mathcal O_{X}$-linear \emph{on both sides}, we may assume that $f = 1$.
We are therefore reduced to checking that $P \bullet 1 =\Phi(P)$ which follows from the definition of the map.
 $\Box$
 
 For example, when $m = 0$, in the case of the affine line with parameter $t$ and derivation $\partial$ and if we choose the usual lifting of Frobenius $t \mapsto t^p$, we have for $k >0$,
$$
\partial \bullet t^{k} = kt^{k-1} - t^k  t^{p-1}\partial^p = (k - t^{p}\partial^p)t^{k-1}.
$$

Thus, if we use $(1, t, \ldots, t^{p-1})$ as a basis for $k[t][[\partial^p]]$ on $k[t^p][[\partial^p]]$, we see that $\partial$ acts as
$$
\left[\begin{array}{ccccc}0 & - t^{p}\partial^p + 1 & 0 & \cdots & 0 \\\vdots & 0 & \ddots & \ddots & \vdots \\\vdots & \vdots & \ddots & \ddots & 0 \\0 & 0 & \cdots & 0 & - t^{p}\partial^p + p-1   \\ - \partial^p & 0 & \cdots & \cdots & 0\end{array}\right].
$$

\textbf{Warning}: Scalar restriction to $\widehat {\mathcal Z}_{X}^{(m)}$ of the action of $\widehat {\mathcal D}_{X}^{(m)}$ on $\widehat{\mathcal Z}\mathcal O^{(m)}_{X}$ is \emph{different} from the natural action of $\widehat {\mathcal Z}_{X}^{(m)}$.

The last result shows that the action of $\widehat {\mathcal D}_{X}^{(m)}$ on $\widehat {\mathcal Z}\mathcal O_{X}^{(m)}$ and the Frobenius of $\widehat {\mathcal D}_{X}^{(m)}$ completely determine each other.
For computations, since the action is $\widehat {\mathcal Z}_{X}^{(m)}$-linear and that $\widehat {\mathcal D}_{X}^{(m)}$ is generated by the operators of order at most $p^m$, the following result might be useful.

\begin{prop}
If $P \in \widehat {\mathcal D}_{X}^{(m)}$ has order at most $p^m$, and $f \in \mathcal O_{X}$, then
$$
P \bullet f  = P(f) + f \Phi(P).
$$
\end{prop}

\textbf{Proof :}
This is a local question and we may therefore assume that we have local coordinates $t_{1}, \ldots, t_{r}$.
By linearity, it is sufficient to show that for $|\underline n| \leq p^m$, we have
$$
\underline \partial^{\langle \underline n \rangle} \bullet f = \underline \partial^{\langle \underline  n \rangle }(f) + f\Phi(\underline  \partial^{\langle \underline  n \rangle }).
$$
We proceed as in the proof of \ref{locdesc}.
By definition, the image of $(1 \otimes f) \otimes (\underline{dt})^{[\underline k]}$ under the morphism
$$
\mathcal O_{X \times_{X'} X} \otimes_{\mathcal O_{X'}}\Gamma_{\bullet}\Omega^1_{X'} \to  \mathcal P_{X,m}.
$$
is exactly
$$
\epsilon(f)\left(\frac 1{p!} \widetilde F^* (\widetilde {\underline \tau'})\right)^{[\underline k]}
$$
where $\epsilon$ denotes the Taylor series.

The coefficient of $\underline\tau^{\{\underline n\}}$ in this series is the same as the coefficient of $\underline \partial^{\{\underline kp^{m+1}\}}$ in $\underline \partial^{\langle \underline n \rangle} \bullet f$.
We may work modulo $\mathcal I^{\{p^m+1\}}$ on the right and we know that, then, all coefficients in $\frac{ \widetilde F^*}{p!}(\widetilde {\underline \tau'})$ are zero unless $\underline n = p^m\mathbf 1_{i}$.
Thus, we are left with the case $\underline k = \underline 0$ that gives $\underline \partial^{\langle \underline  n \rangle }(f)$ and the cases $\underline k = \mathbf 1_{i}$ that gives the different summands of $f\Phi(\underline  \partial^{\langle \underline  n \rangle })$.
$\Box$

Note that the formula is more complicated in higher order: for example if $X$ is the affine line over $\mathbf F_{2}$, $m = 0$ and $\widetilde F$ is the usual lifting of Frobenius, we have
$$
\partial^3 \bullet f = \partial(f) \Phi(\partial^2) + f \Phi(\partial^3).
$$

\section{Higgs modules}\label{higgs}

If $u : X \to Y$ is a morphism of schemes, $\mathcal E$ is an $\mathcal O_{X}$-module and
$$
\theta : \mathcal E \to \mathcal E \otimes_{\mathcal O_{X}} u^*\Omega^1_{Y}
$$
is an $\mathcal O_{X}$-linear map, we will write
$$
\theta^{(1)} : \mathcal E \otimes_{\mathcal O_{X}} u^*\Omega^1_{Y} \stackrel {\theta \otimes \mathrm{Id}}\to \mathcal E \otimes_{\mathcal O_{X}} u^*\Omega^1_{Y} \otimes_{\mathcal O_{X}} u^*\Omega^1_{Y} \stackrel {\mathrm{Id} \otimes \wedge} \to \mathcal E \otimes_{\mathcal O_{X}} u^*\Omega^2_{Y}.
$$

\begin{dfn} Let $u : X \to Y$ be a morphism of schemes.
Then, a \emph{Higgs $u$-module on $X$} is an $\mathcal O_{X}$-module $\mathcal E$, endowed with an $\mathcal O_{X}$-linear map
$$
\theta : \mathcal E \to \mathcal E \otimes_{\mathcal O_{X}} u^*\Omega^1_{Y}
$$
such that $\theta^{(1)} \circ \theta = 0$.
When $u = \mathrm{Id}_{X}$, we say \emph{Higgs module} on $X$.
When $u = F_{X}$ is a (relative iterated) Frobenius morphism, we say Higgs $F$-module on $X$.
\end{dfn}

Higgs $u$-modules form a category with compatible $\mathcal O_{X}$-linear maps as morphisms.
Exactly as modules with integrable connections may be seen as $\mathcal D$-modules, we can interpret the category of Higgs $u$-modules as a category of modules over a suitable ring.

\begin{prop} Let $u : X \to Y$ be a morphism of schemes with $Y$ smooth, and $\mathcal E$ an $\mathcal O_{X}$-module endowed with an $\mathcal O_{X}$-linear map
$$
\theta : \mathcal E \to \mathcal E \otimes_{\mathcal O_{X}} u^*\Omega^1_{Y}.
$$
Then, $\mathcal E$ is a Higgs $u$-module if and only if the dual action
$$
\xymatrix@R=0cm{u^*\mathcal T_{Y} \times \mathcal E \ar[r] & \mathcal E \\ (\xi, s) \ar@{|->}[r] & \xi s}
$$
extends to a structure of $u^*S^\bullet \mathcal T_{Y}$-module.
This is an equivalence of categories.
\end{prop}

\textbf{Proof :}
Note first that a structure of $u^*S^\bullet \mathcal T_{Y}$-module is given by a homomorphism
$$
S^\bullet u^*\mathcal T_{Y} = u^*S^\bullet \mathcal T_{Y} \to \mathcal End_{\mathcal O_{X}}(\mathcal E)
$$
of $\mathcal O_{X}$-algebras.
The universal property of the symmetric algebra tells us that this is equivalent to a linear map
$$
\rho : u^*\mathcal T_{Y} \to \mathcal End_{\mathcal O_{X}}(\mathcal E)
$$
satisfying $\rho(\xi) \circ \rho(\xi') = \rho(\xi') \circ \rho(\xi)$ whenever $\xi, \xi' \in \mathcal T_{Y}$.
Alternatively, it corresponds to a bilinear map
$$
\xymatrix@R=0cm{u^*\mathcal T_{Y} \times \mathcal E \ar[r] & \mathcal E \\ (\xi, s) \ar@{|->}[r] & \xi s}
$$
satisfying $\xi \xi' s = \xi' \xi s$ for $\xi, \xi' \in \mathcal T_{Y}$ and $s \in \mathcal E$.
And the corresponding map
$$
\theta : \mathcal E \to \mathcal E \otimes_{\mathcal O_{X}} u^*\Omega^1_{Y}.
$$
is given in local coordinates by $\theta(s) = \sum \xi_{i} s \otimes \mathrm dx_{i}$.
Now, we compute
$$
\theta^{(1)} (\theta (s)) = \sum_{i<j} (\xi_{i} \xi_{j} s - \xi_{j} \xi_{i} s) \otimes \mathrm dx_{i} \wedge \mathrm dx_{j}.
$$
And we see that this is zero if and only if we always have $\xi_{i} \xi_{j} s = \xi_{j} \xi_{i} s$.
This is exactly the condition we were looking for.
$\Box$

Note that, given a morphism $u : X \to Y$, there is an obvious pull-back morphism $u^*$ from Higgs modules on $Y$ to Higgs $u$-modules on $X$.
There is also a restriction map from Higgs $u$-modules on $X$ to Higgs modules on $X$ that has no interest for us.

\begin{dfn}
Let $X$ be a smooth scheme over a fixed scheme $S$ of characteristic $p$ and $F_{X}$ the relative $m+1$-st iterated Frobenius on $X$.
If $\mathcal E$ is a $\mathcal D^{(m)}_{X}$-module, its \emph{underlying Higgs $F$-module}
is the one obtained by restriction along the $p^m$-curvature map
$$
F_{X}^*S^\bullet\mathcal T_{X'} \to \mathcal D^{(m)}_{X}.
$$
Its \emph{$p^m$-curvature} is the corresponding $\mathcal O_{X}$-linear morphism
$$
\theta : \mathcal E \to \mathcal E \otimes_{\mathcal O_{X}} F_{X}^*\Omega^1_{X'}.
$$
\end{dfn}

Note that this definition of $p^m$-curvature is consistent with definition 3.1.1 of \cite{LeStumQuiros97}.
More precisely, having $p^m$-curvature equal to zero has the same meaning as in \cite{LeStumQuiros97}.

Note that, locally, by definition, we have
$$
\theta(s) = \sum_{i} \partial_{i}^{\{p^{m+1}\}}(s) \otimes \mathrm dt'_{i}.
$$

\begin{dfn}
\begin{enumerate}
\item Let $\mathcal A$ be a ring (in a topos) and $\mathcal I$ a left ideal in $\mathcal A$.
An $\mathcal A$-module $\mathcal E$ is \emph{quasi-nilpotent} if given any section $s \in \mathcal E$, there exists locally $N \in \mathbf N$ such that $\mathcal I^{N} s = 0$.
\item
Let $u : X \to Y$ be a morphism of schemes over some other scheme $S$.
Then, a Higgs $u$-module $\mathcal E$ is \emph{quasi-nilpotent} if it is so as $u^*S^\bullet \mathcal T_{Y}$-module (with respect to its augmentation ideal).
\item
Let $X$ is a smooth scheme over a scheme $S$ of characteristic $p$ and $\mathcal E$ a $\mathcal D^{(m)}_{X}$-module.
Then $\mathcal E$ has \emph{quasi-nilpotent $p^m$-curvature} if the underlying Higgs $F_{X}$-module is quasi-nilpotent.
\end{enumerate}
\end{dfn}

\begin{prop}
Let $X$ be a smooth scheme over $S$ of characteristic $p$ and $\mathcal E$ a $\mathcal D^{(m)}_{X}$-module.
Then, the following are equivalent:
\begin{enumerate}
\item $\mathcal E$ is quasi-nilpotent (with respect to the ideal $\mathcal K_{X}^{(m)}$)
\item $\mathcal E$ is quasi-nilpotent (in the sense of \cite{Berthelot96}, section 2.3)
\item $\mathcal E$ has quasi-nilpotent $p^m$-curvature.
\end{enumerate}
\end{prop}

\textbf{Proof :} All the definitions are local in nature and we may therefore assume that we have local coordinates $t_{1}, \ldots, t_{r}$ on $X$.
Then the first and the third conditions both mean that, locally, we have
$$
(\underline {\partial}^{\langle p^{m+1}\rangle})^{\underline N} (s) = 0 \quad \mathrm{for} \quad |N| >> 0.
$$
Also, the second condition says that, locally again, we have $\partial_{i}^{\langle N \rangle}(s) = 0$ for $N >> 0$.
This is equivalent to the first one because we always have for $t < p^{m+1}$,
$$
\partial_{i}^{\langle kp^{m+1}+t \rangle} = \partial_{i}^t(\partial_{i}^{\langle p^{m+1}\rangle})^k.
\quad \Box$$

Note also that the categories of quasi-nilpotent modules on a ring $\mathcal A$ and on its completion $\widehat {\mathcal A}$ are identical.
In particular, we can always consider the category of quasi-nilpotent $\mathcal D_{X}^{(m)}$-modules as a full subcategory of the category of $\widehat {\mathcal D}_{X}^{(m)}$-module.

In order to go further, we will need the following standard result:

\begin{lem} \label{classeq}
Let $\mathcal R$ be a commutative ring, $\mathcal M$ be a locally free $\mathcal R$ module of finite rank, and $\mathcal A := \mathcal End_{\mathcal R}(\mathcal M)$.
Then, the functors
$$
\mathcal E \mapsto \mathcal Hom_{\mathcal A}(\mathcal M, \mathcal E)
$$
from $\mathcal A$-modules to $\mathcal R$-modules and
$$
\mathcal F \mapsto \mathcal M \otimes_{\mathcal R} \mathcal F
$$
are quasi-inverse to each other.
\end{lem}

\textbf{Proof :} This follows from the fact that the canonical maps
$$
\mathcal M \otimes_{\mathcal R}  \mathcal Hom_{\mathcal A}(\mathcal M, \mathcal E) \to \mathcal E
$$
and
$$
\mathcal F \to \mathcal Hom_{\mathcal A}(\mathcal M, \mathcal M \otimes_{\mathcal R} \mathcal F)
$$
are both bijective.
$\Box$

We denote by $X$ a smooth scheme over a fixed scheme $S$ of characteristic $p$ and let  $F_{X}$ be the relative $m+1$-st iterated Frobenius on $X$.
We also fix a lifting $\widetilde S$ of $S$ modulo $p^2$ as well as a strong lifting $\widetilde F : \widetilde X \to \widetilde X'$ of $F_{X}$ over $\widetilde S$.
Associated to $\widetilde F$, we may consider the Frobenius $\Phi$ of $\widehat {\mathcal D}^{(m)}_{X}$ introduced in the previous section.
Then, we have the following.

\begin{prop} \label{bigequi}
There is an equivalence of categories between $\widehat {\mathcal D}_{X}^{(m)}$-modules and 
$\widehat {S^{\bullet}\mathcal T_{X'}}$-modules given by
$$
\mathcal E \mapsto (F_{X*}\mathcal E)^{1 - \Phi}
$$
and
$$
\mathcal F \mapsto F_{X}^{*}\mathcal F.
$$
\end{prop}

In order to understand this statement, it is necessary to make precise the definition of Frobenius invariants :
$$
(F_{X*}\mathcal E)^{1 - \Phi} := \{ s \in \mathcal E, \quad \forall P \in \widehat {\mathcal D}_{X}^{(m)}, \Phi(P)(s) = P(s)\}.
$$

\textbf{Proof :} Using lemma \ref{classeq}, it follows from Theorem \ref{mainth} that there is an equivalence of categories between $\widehat {\mathcal D}_{X}^{(m)}$-modules and $F_{X*}\widehat {\mathcal Z}_{X}^{(m)}$-modules given by
$$
\mathcal E \mapsto F_{X*}\mathcal Hom_{\widehat {\mathcal D}_{X}^{(m)}}(\widehat {\mathcal Z}\mathcal O_{X}^{(m)}, \mathcal E)
$$
and
$$
\mathcal F \mapsto \widehat {\mathcal Z}\mathcal O_{X}^{(m)}\otimes_{\widehat {\mathcal Z}_{X}^{(m)}} \mathcal F.
$$
We want to identify the right hand sides with the expressions in the proposition.
For the second one, this is easy because
$$
F_{X}^{*} \mathcal F  \simeq F_{X}^{*}F_{X*} \widehat {\mathcal Z}_{X}^{(m)} \otimes_{\widehat {\mathcal Z}_{X}^{(m)}} \mathcal F  \simeq \widehat {\mathcal Z}\mathcal O_{X}^{(m)} \otimes_{\widehat {\mathcal Z}_{X}^{(m)}} \mathcal F.
$$
In order to do the first one, we first check that the canonical map
$$
\xymatrix@R=0cm{\mathcal Hom_{\widehat {\mathcal D}_{X}^{(m)}}(\widehat {\mathcal Z}\mathcal O_{X}^{(m)}, \mathcal E)  \ar[r] & \mathcal E \\ \varphi \ar@{|->}[r] & \varphi(1)}
$$
is injective.
If $\varphi(1) = 0$, then $\varphi(f) = 0$ for any $f \in \mathcal O_{X}$ by linearity.
Moreover, we have the following inclusion
$$
\mathcal Hom_{\widehat {\mathcal D}_{X}^{(m)}}(\widehat {\mathcal Z}\mathcal O_{X}^{(m)}, \mathcal E)  \subset \mathcal Hom_{\widehat {\mathcal Z}_{X}^{(m)}}(\widehat {\mathcal Z}\mathcal O_{X}^{(m)}, \mathcal E)
\simeq \mathcal Hom_{\mathcal O_{X'}}(\mathcal O_{X},\mathcal E)
$$
and it follows that $\varphi =0$.

Now, if $s \in \mathcal E$, we see that the corresponding map
$$
\xymatrix@R=0cm{\varphi : \widehat {\mathcal Z}\mathcal O_{X}^{(m)} \ar[r] &  \mathcal E \\ fQ \ar@{|->}[r] & fQ(s)}
$$
is $\widehat {\mathcal D}_{X}^{(m)}$-linear if and only if $(\Phi(Pf)Q)(s) = P((fQ)(s))$ for all $P \in \widehat {\mathcal D}_{X}^{(m)}$, $f \in \mathcal O_{X}$ and $Q \in \widehat {\mathcal Z}^{(m)}$. 
We may assume that $Q = 1$ and $f = 1$ and this shows that
$$
(F_{X*}\mathcal E)^{1 - \Phi} = F_{X*}\mathcal Hom_{\widehat {\mathcal D}_{X}^{(m)}}(\widehat {\mathcal Z}\mathcal O_{X}^{(m)}, \mathcal E).
\quad \Box$$

\textbf{Warning}: The induced action of $\widehat {\mathcal Z}_{X}^{(m)}$ on $\widehat {\mathcal Z}\mathcal O_{X}^{(m)}$ is \emph{not} the natural one: we need to compose with the automorphism induced by $\Phi$ on $\widehat {\mathcal Z}_{X}^{(m)}$.

We may now state the main theorem of this section (see \cite{OgusVologodsky07}, Theorem 2.8):

\begin{thm} Let $X$ be a smooth scheme over $S$ of positive characteristic $p$.
If there exists a strong lifting of the $m+1$-st iterate of the relative Frobenius of $X$ modulo $p^2$, then there is an equivalence of categories between quasi-nilpotent $ {\mathcal D}_{X}^{(m)}$-modules and quasi-nilpotent Higgs modules on $X'$ given by
$$
\mathcal E \mapsto (F_{X*}\mathcal E)^{1 - \Phi}
$$
and
$$
\mathcal F\mapsto F_{X}^{*} \mathcal F.
$$
\end{thm}

\textbf{Proof :} The functor $F_{X*}$ induces an equivalence of categories between $\mathcal Z_{X}^{(m)}$-modules and  Higgs modules on $X'$.
The theorem therefore follows from Proposition \ref{bigequi} since quasi-nilpotency is clearly preserved under the equivalence.
$\Box$

\textbf{Warning:} Under this equivalence, the $p^m$-curvature of $\mathcal E$ is not obtained as the simple pull-back of the Higgs structure of $\mathcal F$ along Frobenius, we also have to compose with the automorphism induced by $\Phi$ on the center.

Finally, we can give local formulas:

\begin{prop}\label{hgive}
If we are given local coordinates $t_{1}, \ldots, t_{r}$, and $\mathcal F$ is a Higgs module on $X'$ with
$$
\theta(s) = \sum_{i=1}^r \xi'_{i} s \otimes \mathrm dt_{i},
$$
then the $\mathcal D^{(m)}_{X}$-structure of $F_{X}^{*} \mathcal F$ is given by
$$
\partial_{i}^{[p^l]} (1 \otimes s) =
\left\{\begin{array}{l l} 0 & \mathrm{if} \quad l < m
\\
\frac 1{p!} \sum_{j=1}^r \partial_{i}^{[p^m]}(\widetilde F^*(\tilde {t'_{j}})) \otimes \xi'_{j} s & \mathrm{if} \quad l = m.
\end{array}\right.
$$
\end{prop}

\textbf{Proof :}
Follows from proposition \ref{locdesc}.
$\Box$

Unfortunately, it is much more complicated to recover $\mathcal F$ from $\mathcal E$.

\section{Informal complements}\label{comp}

\subsection{Linearizing with respect to the center}\label{example}

As usual, $X$ denotes a smooth scheme over some scheme $S$ of characteristic $p$ and  $F_{X}$ is the relative $m+1$-st iterated Frobenius on $X$.
We also fix a lifting $\widetilde S$ of $S$ modulo $p^2$ as well as a strong lifting $\widetilde F : \widetilde X \to \widetilde X'$ of $F_{X}$ over $\widetilde S$.

First of all, it is important to notice that, in the following commutative diagram, where all maps are $\mathcal O_{X}$-linear,
$$
\xymatrix{F_{X}^*\Gamma_{\bullet}\Omega^1_{X'} \ar[r] \ar@{-->}[d] &  \mathcal P_{X,m} \ar[d] \\F_{X}^*\Gamma_{\bullet}\Omega^1_{X'} \ar[r]^-\simeq &  \mathcal P_{X,m}/\mathcal I\mathcal P_{X,m},}
$$
 the dotted arrow is \emph{not} the identity although both horizontal arrows are defined by the same map $\frac 1{p!} \widetilde F$ on $\Omega^1_{X'}$ and are compatible with divided powers.
This is because the left arrow is not compatible with divided powers.
However, one can show that the dotted arrow is bijective (it is dual to the ring automorphism induced by $\Phi$ on $\widehat{\mathcal Z}\mathcal O^{(m)}_{X}$ - see below).

Recall now that our Frobenius map
$$
\Phi: \widehat {\mathcal D}_{X}^{(m)} \longrightarrow \widehat{\mathcal Z}\mathcal O^{(m)}_{X} \hookrightarrow\widehat {\mathcal D}_{X}^{(m)}
$$
induces the identity on $\mathcal O_{X}$ and a (non trivial) ring automorphism of $\widehat{\mathcal Z}^{(m)}_{X}$.
In particular, $\Phi$ restricts to an automorphism of $\widehat {\mathcal Z}\mathcal O_{X}$ and we will abusively denote by $\Phi^{-1}$ its inverse.
We may then compose $\Phi^{-1}$ on the right with $\Phi$ in order to obtain a modified version of Frobenius 
$$
\widetilde \Phi : \widehat {\mathcal D}_{X}^{(m)} \longrightarrow \widehat{\mathcal Z}\mathcal O^{(m)}_{X} \stackrel {\Phi^{-1}}\simeq \widehat{\mathcal Z}\mathcal O^{(m)}_{X} \hookrightarrow\widehat {\mathcal D}_{X}^{(m)}.
$$
that will induce the identity on $\widehat{\mathcal Z}\mathcal O^{(m)}_{X}$.

Using this modified Frobenius gives a twisted version of our fundamental algebra isomorphism
$$
\widehat {\mathcal D}_{X}^{(m)} \simeq \mathcal End_{\widehat {\mathcal Z}_{X}^{(m)}}(\widehat{\mathcal Z}\mathcal O^{(m)}_{X})
$$
which is now the identity on the center $\widehat {\mathcal Z}_{X}^{(m)}$.
The corresponding $\widehat {\mathcal D}_{X}^{(m)}$-module structure is given by a similar formula
$$
P \widetilde \bullet fQ = \widetilde \Phi (Pf)Q.
$$
This twisted action is given by more complicated formulas but has the advantage to be a true Azumaya splitting for $\widehat {\mathcal D}_{X}^{(m)}$.
Of course, one can use this twisted action in order to obtain an equivalence between quasi-nilpotent $\mathcal D^{(m)}_{X}$-modules and Higgs modules on $X$.
Actually, $\Phi$ induces an autoequivalence of the category of quasi-nilpotent Higgs modules on $X'$ and the twisted equivalence is obtained by composition with the old one.

\subsection{van der Put's construction}

Marius van der Put only deals with differential operators of level zero and we will therefore stick here to this case and drop $m$ from the notations.
He actually works in the theory of differential fields (\cite{Put95}, \cite{PutSinger97}) but the translation into our language is straightforward.
Also, we will not follow his notations but only try to explain his clever approach.
We let $X$ be a curve over a perfect field, with a fixed coordinate $t$ and denote as usual the corresponding derivation by $\partial$.

It is important to note that van der Put works with the twisted version of the theory (see section \ref{example}).
The first point is to remark that the Azumaya nature of $\widehat {\mathcal D}_{X}$ is completely described by the element
$$
H := \widetilde \Phi (\partial) \in \widehat{\mathcal Z}\mathcal O_{X}.
$$
In what follows, we make $\widehat {\mathcal D}_{X}$ act on $\widehat{\mathcal Z}\mathcal O_{X}$ by extending trivially the action on $\mathcal O_{X}$ (via the isomorphism $\mathcal O_{X}\otimes_{\mathcal O_{X'}} \widehat{\mathcal Z}_{X} \simeq \widehat{\mathcal Z}\mathcal O_{X}$).
Now, since the twisted action is the identity on $\widehat{\mathcal Z}_{X}$, we must have $\partial^p \widetilde \bullet 1 = \partial^p$ and this leads to the condition
$$
\partial^{p-1}(H) + H^p = \partial^p.
$$
Thus, van der Put's question reduces to finding such an element.
He does it by successive approximations, lifting an action on the quotient $\mathcal O_{X}$ (see \cite{Put95}, Lemma 1.6).
The point is to notice that such an action corresponds to an element $h \in \mathcal O_{X}$ such that $\partial^{p-1}(h) + h^p = 0$.
In other words, he lifts the quotient isomorphism :
$$
\xymatrix{\widehat {\mathcal D}_{X} \ar[rr]^-\simeq \ar@{->>}[d] &&\mathcal End_{\widehat {\mathcal Z}_{X}}(\widehat{\mathcal Z}\mathcal O_{X}) \ar@{->>}[d] \\ \mathcal D_{X}/\mathcal K_{X} \ar[rr]^-\simeq && \mathcal End_{\mathcal O_{X'}}(\mathcal O_{X}).
}
$$
We will now check that he obtains exactly the map induced by $\widetilde \Phi$ if we can use $t \mapsto t^p$ as lifting of Frobenius modulo $p^2$. First of all, we have
$$
\Phi(\partial) = - t^{p-1}\partial^p \quad \mathrm{and} \quad  \Phi(\partial^p) = \partial^p - t^{p(p-1)}\partial^{p^2}.
$$
Since $\Phi$ induces a morphism of rings on $\widehat{\mathcal Z}\mathcal O_{X}$, we also have for all $k \geq 0$
$$
\Phi(\partial^{p^k}) = \partial^{p^k} -(t^{p-1})^{p^k}\partial^{p^{k+1}}
$$
from which we derive
$$
\Phi^{-1}(\partial^p) = \sum_{k=1}^{\infty} (t^{p-1})^{\frac {p^{k} -1}{p-1} -1} \partial^{p^k} = \sum_{k=1}^{\infty} t^{p(p^{k-1} -1)} \partial^{p^k}.
$$
Since $\Phi$ is $\mathcal O_{X}$-linear, we can deduce
$$
H := \widetilde{\Phi}(\partial) =  - t^{p-1} \Phi^{-1}(\partial^p) =  - \sum_{k=1}^{\infty} t^{(p^{k} -1)} \partial^{p^k}.
$$
This is exactly the same (with different notations) as Example 1.6.1 of \cite{Put95}.
The reader can also look at the proof of theorem 13.5 of \cite{PutSinger97}).
Note also that many examples are worked out in this book.

\subsection{Working before completion}

We keep the same notations as in section \ref{example}.
It is not difficult to see that the map $\Phi$ exists before completion.
One just has to be careful when dualizing.
This is possible because we have explicit local formulas.
Thus, there is a well defined map
$$
\Phi: \mathcal D_{X}^{(m)} \longrightarrow \mathcal Z\mathcal O^{(m)}_{X} \hookrightarrow\mathcal D_{X}^{(m)}.
$$
Using our local formulas again, it is not difficult either to see that $\Phi$ sends $\mathcal Z_{X}$ inside itself.
But the induced map need not be bijective (in particular, we cannot define $\widetilde \Phi$ before completion).

For the same reasons, there exists a morphism of $\mathcal O_{X}$-algebras at finite level
$$
\mathcal D_{X}^{(m)} \to \mathcal End_{\mathcal Z_{X}^{(m)}}(\mathcal Z\mathcal O^{(m)}_{X}).
$$
But again, this is not an isomorphism.
However, this map is semi-linear with respect to the endomorphism $\Phi$ of $\mathcal Z_{X}$ and linearization provides an isomorphism
$$
\mathcal Z_{{X}_{\nwarrow \atop \Phi}}\!\!\! \otimes_{\mathcal Z_{X}} \mathcal D_{X}^{(m)} \simeq \mathcal End_{\mathcal Z_{X}^{(m)}}(\mathcal Z\mathcal O^{(m)}_{X}).
$$
Note that the fact that this is an isomorphism follows from our previous results because it is a morphism of locally free $\mathcal Z_{X}$-modules and that it is therefore sufficient to prove bijectivity after completing.

The endomorphism induced by $\Phi$ on $\mathcal Z_{X} \simeq S^\bullet \mathcal T_{X'}$ corresponds to a morphism
$$
\alpha : \check{\mathbf T}_{X'} \to \check{\mathbf T}_{X'}.
$$
We did not check it in general, but it is likely that $\alpha$ is surjective \'etale and provides an \emph{\'etale} Azumaya splitting of $\mathcal D_{X}^{(m)}$ (see Proposition 2.5 (1) of \cite{OgusVologodsky07}  for the case $m = 0$).

Finally, we can pullback Higgs modules along Frobenius in order to get a $\mathcal D_{X}^{(m)}$-module, but it becomes an equivalence only when restricting to quasi-nilpotent objects.

\subsection{Glueing}

As before, $X$ denotes a smooth scheme over some scheme $S$ of characteristic $p$ and  $F_{X}$ is the relative $m+1$-st iterated Frobenius on $X$.
We also fix a lifting $\widetilde S$ of $S$ modulo $p^2$ as well as a lifting $\widetilde X/\widetilde S$ of $X/S$.

Assume for the moment that we are given two strong liftings $\widetilde F_{1}, \widetilde F_{2} : \widetilde X \to \widetilde X'$ of $F_{X}$ over $\widetilde S$ (we may always twist a lifting of $X$ and therefore assume that $\widetilde F_{1}$ and $\widetilde F_{2}$ have the same domain).
Considering local descriptions again, one sees that the map
$$
\widetilde F_{2} - \widetilde F_{1} : \mathcal O_{\widetilde X'} \to \mathcal O_{\widetilde X}
$$
induces a derivation
$$
\frac 1 {p!} (\widetilde F_{2} -\widetilde F_{1}) :\mathcal O_{X'} \to \mathcal O_{X}.
$$
and we obtain a natural $\mathcal O_{X'}$-linear map
$$
u_{12} = \frac 1 {p!} (\widetilde F_{2} -\widetilde F_{1}) : \Omega^1_{X'} \to \mathcal O_{X}.
$$
Moreover, if $\widetilde F_{3}$ is another strong lifting of $F_{X}$, we have (with the obvious analogous notation) $u_{13} = u_{12} + u_{23}$.

The map $u_{12}$ extends to a morphism of $\mathcal O_{X'}$-algebras $S^\bullet \Omega^1_{X'} \to \mathcal O_{X}$ and by linearity to map $F^*_{X}S^\bullet \Omega^1_{X'} \to \mathcal O_{X}$ and we finally embed into $F^*_{X}S^\bullet \Omega^1_{X'}$ in order to get an endomorphism of the $\mathcal O_{X}$-algebra $F^*_{X}S^\bullet \Omega^1_{X'}$ that we still call $u_{12}$.
We may then define $\phi_{12} = \mathrm {Id} - u_{12}$ and remark that we have $\phi_{12}\phi_{23} = \phi_{13}$.
It follows that, when there are only local liftings of Frobenius, we obtain glueing data and we will call $\mathcal A^{(m)}_{X}$ the resulting $\mathcal O_{X}$-module which is locally isomorphic to $F^*_{X}S^\bullet \Omega^1_{X'}$.

In order to dualize, we introduce the completed divided power envelopes
$$
\widehat {\mathcal Z^{(m)\mathrm{DP}}_{X}}, \quad \widehat {\mathcal Z\mathcal O^{(m)\mathrm{DP}}_{X}} \quad \mathrm{and} \quad \widehat {\mathcal D_{X}^{(m)DP}}
$$
of
$$
\mathcal Z^{(m)}_{X}, \quad \mathcal Z \mathcal O^{(m)}_{X} \quad \mathrm{and} \quad \mathcal D_{X}^{(m)}
$$
respectively along the augmentation ideal of $\mathcal Z^{(m)}_{X}$.
Note that we have
$$
\widehat {\mathcal Z^{(m)\mathrm{DP}}_{X}} \simeq \widehat {\Gamma}_{\bullet}\mathcal T_{X'}, \quad 
\widehat {\mathcal Z\mathcal O^{(m)\mathrm{DP}}_{X}} \simeq \widehat {\mathcal Z^{(m)\mathrm{DP}}_{X}}\otimes_{\mathcal Z^{(m)}_{X}} \mathcal Z \mathcal O_{X}^{(m)}
$$
and
$$
\quad \widehat {\mathcal D_{X}^{(m)\mathrm{DP}}} \simeq \widehat {\mathcal Z^{(m)\mathrm{DP}}_{X}}\otimes_{\mathcal Z^{(m)}_{X}} \mathcal D^{(m)}_{X}.
$$
Thus, dual to $\mathcal A^{(m)}_{X}$ we obtain an $\mathcal O_{X}$-module $\mathcal B^{(m)}_{X}$, which is in fact a locally free $\widehat {\mathcal Z\mathcal O^{(m)\mathrm{DP}}_{X}}$-module of rank $1$.
In order to verify this, it is sufficient to check that, locally, the dual glueing map
$$
\check \phi_{12} : \widehat {\mathcal Z\mathcal O^{(m)\mathrm{DP}}_{X}} \to \widehat {\mathcal Z\mathcal O^{(m)\mathrm{DP}}_{X}}
$$
is $\widehat {\mathcal Z\mathcal O^{(m)\mathrm{DP}}_{X}}$-linear:
it is $\mathcal O_{X}$-linear by definition;
moreover, by construction, $\check \phi_{12}$ acts as the identity on the augmentation ideal of $\widehat {\mathcal Z^{(m)\mathrm{DP}}_{X}}$ because $u_{12}$ sends the augmentation ideal of $S^\bullet \Omega^1_{X'}$ inside $\mathcal O_{X}$.

Note that, locally, $\widetilde F_{2} -\widetilde F_{1}$ sends $\widetilde{\mathcal I}$ into $p\mathcal O_{\widetilde X \times_{\widetilde X'} \widetilde X}$ and that we obtain a map
$$
\xymatrix@R0cm{\Omega^1_{X'} \ar[r]^-{U_{12}} & \mathcal O_{X \times_{X'} X} \\ \omega' \ar@{|->}[r] & u_{12}(\omega') \otimes 1 - 1 \otimes u_{12}(\omega').}
$$
A construction analog to the above provides an automorphism $\mathrm{Id} - U_{12}$ of $\mathcal O_{X \times_{X'} X} \otimes S^\bullet \Omega^1_{X'}$ giving rise to a commutative diagram
$$
\xymatrix{\mathcal O_{X \times_{X'} X} \otimes S^\bullet \Omega^1_{X'} \ar[rd]_{\frac1{p!}{\widetilde F_{1}}} \ar[rr]^\simeq && \mathcal O_{X \times_{X'} X} \otimes S^\bullet \Omega^1_{X'} \ar[ld]^{\frac1{p!}{\widetilde F_{2}}} \\ & \mathcal P_{Xm}&}.
$$
Dualizing provides a morphism
$$
\widehat {\mathcal D}_{X}^{(m)} \to \mathcal End_{\widehat {\mathcal Z_{X}^{(m)\mathrm{DP}}}}(\mathcal B^{(m)}_{X}),
$$
and extending scalars, a map which is an isomorphism
$$
\widehat {\mathcal D_{X}^{(m)\mathrm{DP}}} \simeq \mathcal End_{\widehat {\mathcal Z_{X}^{(m)\mathrm{DP}}}}(\mathcal B^{(m)}_{X}),
$$
as one easily checks.

\subsection{Ogus-Vologodsky's construction}

We keep the same notations as before but we assume that $m = 0$ and we drop it from the notations.
Also, in what follows, there will be some sign differences with the article of Ogus and Vologodsky (\cite{OgusVologodsky07}), due to the fact that they choose to divide by $p$ instead of $p!$ and that $p! = -p \mod p^2$.

They start with a global definition of $\mathcal A_{X}$ that we do not want to recall here (\cite{OgusVologodsky07}, Theorem 1.1).
In order to obtain a local description of this object, they introduce the notion of ``splitting of Cartier'' $\zeta : \Omega^1_{X'} \to \Omega^1_{X}$.
When $\widetilde F$ is a lifting of Frobenius, the standard example of such is splitting is given by Mazur's construction and is simply the composite of the map
$$
\frac {\widetilde F} {p!} : \Omega^1_{X'} \to \mathcal I_{X}
$$
with the canonical projection $\mathcal I_{X} \to \mathcal I_{X}/\mathcal I_{X}^{[2]} = \Omega^1_{X}$.
They linearize and dualize in order to obtain a morphism $\phi : \mathcal T_{X} \to F_{X}^*\mathcal T_{X'}$.
Thus, by construction, there is a commutative diagram
$$
\xymatrix{\Phi : \mathcal D_{X} \ar[rr]^\Phi && \mathcal D_{X} \\ \mathcal T_{X} \ar@{^{(}->}[u] \ar[rr]^\phi && F_{X}^*\mathcal T_{X'}\ar@{^{(}->}[u]  }
$$
where the map on the right is the $p$-curvature map.
From $\phi$, they construct a morphism $h : \check{\mathbf T}_{X'} \to \check{\mathbf T}_{X'}$ that is given on local sections by $h^*(\partial_{i}^p) = \Phi(\partial_{i})^p$ if we use as usual the $p$-curvature map to embed $\mathcal T_{X'}$ into $\mathcal D_{X}$.
When $m = 0$, our formula in Corollary \ref{phibar} reads $\Phi(\partial_{i}^p) = \partial_{i}^p + \Phi(\partial_{i})^p$ and it follows that our map $\alpha$ is the same as the one in section 2.2 of \cite{OgusVologodsky07}.
In particular, our theory is fully compatible with theirs even if the approach might sound different.

For example, their formula (2.11.2) reads locally
$$
\nabla(1 \otimes s) = \sum_{i=0}^r \xi_{i}' s \otimes \zeta(\mathrm dt'_{i}).
$$
and we have, with Cartier-Mazur's splitting corresponding to a lifting $\widetilde t'_{i} \mapsto \widetilde t_{i}^p + p \widetilde g_{i}$,
$$
\zeta(\mathrm dt'_{i}) = - t_{i}^{(p-1)} \mathrm dt_{i} - \sum_{j=1}^r \partial_{j}(g_{i}) \mathrm dt_{j}.
$$
It follows that
$$
\partial_{i}(1 \otimes s) = - t_{i}^{(p-1)} \otimes \xi_{i}' s  - \sum_{j=1}^r \partial_{i}(g_{j}) \otimes \xi'_{j} s
$$
which is exactly what we found in Proposition \ref{hgive}.

\subsection{Frobenius descent}
 
When $s \in \mathbf N$, any divided power structure of level $m$ on an ideal $\mathcal I$ is \emph{a fortiori} a divided power structure of higher level $m + s$.
It follows that if $X \hookrightarrow Y$ is any embedding, there exists natural map $P_{Xm}(Y) \to P_{Xm+s}(Y)$.
And from this, one derives, for a smooth scheme $X$ over our fixed scheme $S$, a morphism of rings $\mathcal D_{X}^{(m)} \to \mathcal D_{X}^{(m+s)}.$
It is given on local generators by $\partial_{i}^{[p^l]} \mapsto \partial_{i}^{[p^l]}$ for $l \leq m$ but this is not an injective map in general.
For example, if $S$ has positive characteristic $p$ and $s > 0$, then $\partial_{i}^{\{p^{m+1}\}}$ is sent to zero.

Let us be more specific in this positive characteristic situation.
Given any $s$, we will denote by $X^{(s)}$ the pull-back of $X$ along the $s$-st iteration of Frobenius.
For $s > 0$, the kernel of the map $\mathcal D_{X}^{(m)} \to \mathcal D_{X}^{(m+s)}$ is nothing but the ideal $\mathcal K^{(m)}_{X}$ generated by the image of $\mathcal T_{X^{(m+1)}}$ under the $p^m$-curvature map.
Moreover, it induces an isomorphism of rings
$$
\mathcal End_{\mathcal O_{X^{(m+1)}}}(\mathcal O_{X}) \simeq \mathcal D_{X,p^{m+1}-1}^{(m+s)}
$$
where the right hand side is the subring of differential operators of order less than $p^{m+1}$.

We now come to the basic fact on Frobenius descent : the $s$-th iteration of Frobenius induces a morphism
$$
F^{s*} :P_{Xm+s} \to  P_{X^{(s)}m}.
$$
It is given locally by $\tau_{i}^{\{k\}_{m}} \mapsto \tau_{i}^{\{kp^s\}_{m+s}}$ and we obtain by duality, a morphism
$$
\mathcal D_{X}^{(m+s)} \to F^{s*}\mathcal D_{X^{(s)}}^{(m)}
$$
given locally by $\partial_{i}^{\langle k \rangle} \mapsto 1 \otimes \partial_{i}^{\langle k/p^s \rangle}$ if $p^s | k$ and $0$ otherwise.
Note that $F^{s*}\mathcal D_{X^{(s)}}^{(m)}$ has no natural ring structure so that this morphism cannot be seen as a  morphism of rings.
However, it induces an isomorphism of rings $ \mathcal Z \mathcal O_{X}^{(m+s))} \simeq F^{s*}\mathcal Z \mathcal O_{X^{(s)}}^{(m)}$ compatible with the $p$-curvature maps of level $m+s$ and $m$ respectively.

If we choose compatible strong liftings, there is an obvious commutative diagram
$$
\xymatrix{ && \mathcal P_{X^{(s)}m} \ar[dd]^{F^{s*}}\\ \Omega^1_{X^{(m+s+1)}} \ar[rru]^{\frac 1{p!}\widetilde F^{m+1*}} \ar[rrd]_{\frac 1{p!}\widetilde F^{m+s+1*}}&&\\ && \mathcal P_{Xm+s} }
$$
that gives by duality a commutative diagram (we ignore completions)
$$
\xymatrix{\mathcal D^{(m+s)}_{X} \ar[r] \ar[d] \ar@/^1cm/[rr]^{\Phi^{(m+s)}_{X}} & \mathcal Z \mathcal O_{X}^{(m+s))} \ar[d]^\simeq \ar@{^{(}->}[r] & \mathcal D^{(m+s)}_{X} \ar[d] \\ F^{s*}\mathcal D^{(m)}_{X^{(s)}} \ar[r] \ar@/_1cm/[rr]_{F^{s*}\Phi^{(m)}_{X^{(s)}}} & F^{s*}\mathcal Z \mathcal O_{X^{(s)}}^{(m)} \ar@{^{(}->}[r]  & F^{s*}\mathcal D^{(m)}_{X^{(s)}} .}
$$
For the same reason, we have a commutative diagram
$$
\xymatrix{\mathcal D^{(m+s)}_{X} \ar[r] \ar[d] & \mathcal End_{\mathcal Z_{X}^{(m+s)}}(\mathcal Z\mathcal O_{X}^{(m+s)}) \ar[d] \\ F^{s*}\mathcal D^{(m)}_{X^{(s)}} \ar[r] & F^{s*}\mathcal End_{\mathcal Z_{X^{(s)}}^{(m)}}(\mathcal Z\mathcal O_{X^{(s)}}^{(m)}) .}
$$

Finally, when restricting to quasi-nilpotent objects everywhere, we have compatible equivalences of categories induced by Frobenius pull backs
$$
\xymatrix{
& \textrm {\{Higgs modules on $X^{(m+s+1)}$\}}  \ar[rd]^\simeq_{F^{m+s+1*}}  \ar[ld]^{F^{m+1*}}_\simeq & \\
\textrm{\{${\mathcal D}_{X^{(s)}}^{(m)}$-modules\}} \ar[rr]^\simeq_{F^{s*}} && \textrm{\{${\mathcal D}_{X}^{(m+s)}$-modules\}} }
$$
In particular, we recover our equivalence of categories as the composition of Frobenius descent and the usual case $m = 0$ (again, we restrict to quasi-nilpotent objects) :
$$
\xymatrix{\textrm {\{Higgs modules on $X^{(m+1)}$\}}  \ar[r]^-\simeq_-{F^*} & \textrm{\{${\mathcal D}_{X^{(m)}}^{(0)}$-modules\}} \ar[r]^-\simeq_-{F^{m*}} &\textrm{\{${\mathcal D}_{X}^{(m)}$-modules\}}.}
$$
Note however that, as far as we understand, the Azumaya nature of ${\mathcal D}_{X}^{(m)}$ cannot be derived from the case $m = 0$.

\bibliographystyle{plain}
\bibliography{BiblioBLS}

\bigskip
{\small 
\begin{tabular}{lclcl}
Michel Gros && Bernard Le Stum && Adolfo Quir\'os \\
IRMAR && IRMAR && Departamento de Matem\'aticas \\
Universit\'e de Rennes 1 && Universit\'e de Rennes 1&& Universidsad Aut\'onoma de Madrid\\
Campus de Beaulieu && Campus de Beaulieu&& Ciudad Universitaria de Cantoblanco \\
35 042 Rennes Cedex && 35 042 Rennes Cedex && E-28049 Madrid \\
France && France && Espa\~na
\end{tabular}
}
\end{document}